\documentclass[a4paper,11pt]{article}
\usepackage[utf8]{inputenc}
\usepackage[margin=3cm]{geometry}
\usepackage{amsfonts,amsmath,amssymb,amscd,graphicx}
\usepackage[colorlinks]{hyperref}
\usepackage{enumitem}

\newtheorem{theorem}{Theorem}[section]
\newtheorem{corollary}[theorem]{Corollary}
\newtheorem{lemma}[theorem]{Lemma}

\newtheorem{proposition}[theorem]{Proposition}

\newtheorem{definition}[theorem]{Definition}
\newtheorem{remark}[theorem]{Remark}
\newtheorem{example}[theorem]{Example}
\numberwithin{equation}{section}

\newenvironment{preuve}[1][]
{\vskip 2mm  \emph{\bf Proof#1. }}{$\Box$ \vskip 2mm}

\newcommand{\sn}{\sqrt n \mathbb S^{n-1}}

\newcommand{\Nn}{\mathbb{N}}
\newcommand{\R}{\mathbb{R}}

\newcommand{\ep}{\epsilon}
\newcommand{\si}{\sigma}

\let\epsilon=\varepsilon
\let\ln=\log

\newcommand{\vol}{\ensuremath \mathrm{vol}}
\newcommand{\Tr}{\ensuremath \mathrm{tr}}
\newcommand{\Cone}{\ensuremath \mathrm{Cone}}
\newcommand{\rank}{\ensuremath \mathrm{rank}}
\newcommand{\crit}{\ensuremath \mathrm{Crit}}
\newcommand{\Ind}{\ensuremath \mathrm{ind}}
\newcommand{\sym}{\ensuremath \mathrm{Sym}}

\begin{document}

\title{\bf Asymptotic topology  \\ of excursion and nodal sets \\of  Gaussian random fields}
\author{\sc Damien Gayet}
\maketitle

\begin{abstract}
	Let $M $ be a compact smooth manifold of dimension $n$ with or without boundary, and $ f : M\to \R$ be a smooth Gaussian random field. It is very natural to suppose that for a large positive real $u$, the random excursion set $\{f\geq u\}$  is mostly composed of a union of disjoint topological $n$-balls. Using the constructive part of (stratified) Morse theory we prove that in average, this intuition is true, and provide  for large $u$ the asymptotic of the expected number of such balls, and so of connected components of $\{f\geq u\}$, see Theorem~\ref{sph0}. We similarly show that in average, the high nodal sets $\{f=u\}$ are mostly composed of spheres, with the same asymptotic than the one for excursion set. A refinement of these results using the average of the Euler characteristic given by~\cite{adler} provides a striking asymptotic of the constant defined by F. Nazarov and M. Sodin, again for large $u$, see Theorem~\ref{ns2}. This new Morse theoretical approach of random topology also applies to spherical spin glasses with large dimension, see Theorem~\ref{spin}.
\end{abstract}

Keywords: {Random topology, excursion set, smooth Gaussian field, Morse theory, spin glasses.}

\textsc{Mathematics subject classification  2010}: 60K35, 26E05.


\tableofcontents
\section{Introduction}

\subsection{The results}
\paragraph{Setting and notations.}
Let $M$ be a compact
smooth manifold with or without boundary, or more generally a compact \emph{Whitney stratified set}, a family of sets which contains, for instance, manifolds with corners as affine hypercubes, see Definition~\ref{stratified} below. 
Let $ f :  M\to \R$ be a random centered smooth Gaussian field with constant variance. For any $u\in \R$, denote by 
$ E_{u}(M,f)$ the \emph{excursion set} of $f$ over the threshold $u$, 
or $E_u$ when $f$ is implicit, that is 
$$ E_{u}(M,f)= \{x\in M, f(x) \geq u\}.$$
The \emph{sojourn set} under $u$ is the sublevel
$$ S_{u}(M,f)= \{x\in M, f(x) \leq u\},$$
and the \emph{nodal set} at $u$ is the level set 
$$ Z_u (M,f)= \{x\in M, f(x) = u\}.$$ 
The statistical geometric and topological features of these random sets
have been studied since the 50's, see paragraph~\ref{history} for a survey of past results in this topic. Topological observables of interest are the Euler characteristic of the excursion or level set, its number of connected components, its Betti numbers, or more precisely, the homeomorphic type of its components. The first is local, as the volume, hence has been studied first. The other ones are global, hence more difficult to access, and has been studied more recently.

Classical  Morse theory allows to understand partially the topology of a differential manifold $M$ through the critical points of some unique generic function, see Section~\ref{morsem}. Quite surprisingly, it allows to compute the Euler characteristic through a similar Euler-Morse characteristic involving only critical points, see~(\ref{chi}). This beautiful equality has been extensively used in the probabilistic litterature in order to compute the the average of $\chi(E_u(M,f))$. Morse theory also provides informations about Betti numbers through so-called Morse inequalities, (see Theorem~\ref{milnor} assertion~\ref{weak}). It has also been used to bound above the mean Betti numbers of $Z_u(M,f)$, see \S~\ref{history}. In this paper, we apply another part of Morse theory, which allows to be far more precise, namely to describe the average topological type of the excursion and nodal sets of random functions for large positive levels $u$. Note that for instance, the Euler caracteristic of an circle (and hence of an annulus or a full torus) vanishes, as it is the case for any oriented closed manifold with odd dimension. 

\paragraph{The main results.}
It is very natural to believe that for large positive $u$, most of connected components of the excursion set $E_u$ are emerging islands, sometimes called \emph{bumps} or \emph{blobs} in the litterature, that is components diffeomorphic to the standard $n$-ball, if $M$ has dimension $n$. 
%
In this paper, we prove that this intuition is correct in a quantitatively way. In order to make the statement more formal, 
we follow~\cite{GWBetti}: for any smooth compact smooth submanifold $\Sigma\subset \R^n$, possibly with boundary, and any subset $E$ of a manifold $M$, we set 
\begin{align*}
N_{\Sigma}(E)&=\#\{\text{connected components } B \text { of } E \ | \ B \text{ is diffeomorphic to }\Sigma\}\\
\text{ and } N(E)&=b_0(E)=\#\{\text{connected components of } E\}.
\end{align*}
We emphasize that in the case where $M$ is a manifold with boundary,  we count for $N(M)$ the components of $M$ touching the boundary as well. 

\begin{figure}
	\centering
	\includegraphics[height=0.6\textwidth]{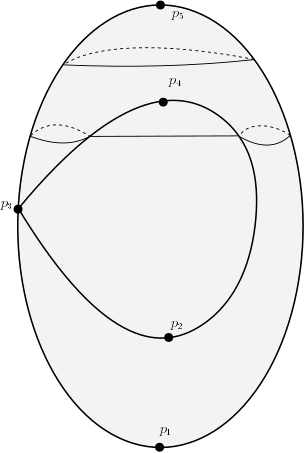}
	\caption{A stratified (pinched) mirror $M\subset \widetilde M\subset \R^3$~\cite[p. 6]{goresky}. It is the union of a pinched torus (the frame of the mirror) and the vertical disc. Here, $p_3$ is $\partial_0M$, the boundary of the vertical mirror without $p_3$ is $\partial_1 M$ and $\partial_2 M$ is the rest, that is the union of the frame without $\partial_0 M$ and $\partial_1 M$ and the open vertical mirror. In this case $M$ is a regular cone space which is not locally convex, see Definition~\ref{locconv}.  The singular point $p_3$ is critical for any function, and the height is a Morse function $f$ for $M$ in the sense of Definition~\ref{defmorse}, for which the $p_i$'s are critical. Here the indices~(\ref{index}) are  $\Ind(p_1)=0$, $\Ind(p_2)=0$ (since the stratum of $p_2$ is $\partial_1 M$, $\Ind(p_4)=1$ and $\Ind(p_5)=2$. Two level lines of $f$ are shown. Stratified Morse theory describes the changes of topology of the sublevel sets of $f$, see Theorem~\ref{Tgoresky}. }\label{tore-bizarre}
\end{figure}

We begin with a corollary. Let $M$ be a compact smooth manifold with or without boundary and $f: M\to \R$ be a centered Gaussian field
satisfying 
the hypotheses (\ref{C1}) (regularity) and (\ref{C2}) (non-degeneracity) given below. Then,  
$f$ induces a metric over $M$ 
by~\cite[(12.2.1)]{adler}:
\begin{equation}\label{metric0}
\forall x\in M, X,Y\in T_x M, \ 
g(X,Y) =\mathbb E(df(x)X, df(x)Y),
\end{equation}
where $df(x)$ denotes the differential of $f$ at $x$.
\begin{corollary}\label{coroman}Let $M$ be a compact $C^3$ manifold with or without boundary, and $f: M\to \R$ be a random centered Gaussian field satisfying conditions (\ref{C1}) (regularity), (\ref{C2})(non-degeneraticity) and (\ref{C4}) (constant variance), and $g$ be de metric defined by~(\ref{metric0}). Then
	$$  \forall u\in \R, \ 
	\mathbb E N_{\mathbb B^n}(E_u(M,f)) =
	\frac{1}{\sqrt{2\pi}^{n+1}} 	\vol_g (M)u^{n-1}e^{-\frac{1}2u^2} 
	\left(1 
	+O_{u\to+\infty}(\frac1u)\right).
	$$
	Here, the error term depends only on the 4-jet of the covariance on the diagonal $M\times M$.
	The same holds for 
	$N(E_u)$, $N_{\mathbb S^{n-1}}(Z_u)$ and $N(Z_u)$ instead of $N_{\mathbb B^n}(E_u)$. 
\end{corollary}
Note that this is the first asymptotic for the average number of components of given diffeomorphism type of random smooth subsets, and the first asymptotic for the number of components in dimension larger than 2. This corollary is a particular case of a far more general theorem, which holds for \emph{Whitney stratified sets}, see Definition~\ref{stratified}:
\begin{theorem}\label{sph0}
	Let $n$ be a positive integer, $\widetilde M$ be a $C^3$ manifold of dimension $n$,  $M\subset \widetilde M$ be a compact $C^2$  Whitney stratified set  of dimension $n$ satisfying conditions~(\ref{hausstrat}) (gentle boundaries) and~(\ref{bettistrat0}) (mild local connectivity). Let $ \widetilde f: \widetilde M \to \R$ be a random centered Gaussian field satisfying conditions (\ref{C1}) (regularity), (\ref{C2}) (non-degeneracity) and (\ref{C4}) (constant variance), $f=\widetilde f_{|M}$ and $g$ be the metric induced by $\tilde f$ and defined by~(\ref{metric}). Then
	\begin{equation}\label{sky}  
\forall u\in \R, \ 	 \mathbb E N_{\mathbb B^n}(E_u(M,f)) =
\frac{1}{\sqrt{2\pi}^{n+1}} 	\vol_g (\partial_n M)u^{n-1}e^{-\frac{1}2u^2} \left(1 
	+O_{u\to+\infty}(\frac1u)\right).
	\end{equation}
		The same holds for 
	$N(E_u)$ instead of $N_{\mathbb B^n}(E_u)$. 
	Here, $\partial_nM$ denotes the stratum of maximal dimension, see~(\ref{dimension}) 	and the error term depends only on the 4-jet of the covariance on the diagonal $\widetilde M\times \widetilde M$.
	
	If moreover $M$ satisfies the further condition~(\ref{zetastrat}) (milder topology), then~(\ref{sky}) holds with $N_{\mathbb S^{n-1}}(Z_u)$ and $N(Z_u)$ instead of $N_{\mathbb B^n}(E_u)$. 
\end{theorem}
Since they need a lot of material, we postpone the necessary definitions and conditions to Section~\ref{morsestrat} and~\ref{proof}. However, let us say here that stratified sets are decomposed into submanifolds which are called \emph{strata} of different dimensions denoted by $\partial_j M$, where $j$ is the dimension.
\begin{example}\label{exmo}
	Manifolds with or without boundaries and affine hypercubes satisfy the hypotheses of Theorem~\ref{sph0}.
	For a manifold without boundary, $\partial_n M=M$ and for any $j\leq n-1$, $\partial_j M= \emptyset$.
	For a manifold with boundary, $\partial_n M=M\setminus \partial M$, $\partial_{n-1}M=\partial M$ and there is no other strata. For the hypercube $[0,1]^n$, $\partial_j M$ is the union of the faces of dimension $j$, and $\partial_n M= \overset{\circ}{M}$. A more exotic example is provided by Figure~\ref{tore-bizarre}.
\end{example}

\begin{remark}\label{rem}
	\begin{enumerate}
		\item A version of Theorem~\ref{sph0} with a precise error bound and for spaces with positive codimension, that is $\dim M<\dim \widetilde M$, is given by Theorem~\ref{sph03}.
		\item Condition~(\ref{C4}) could be dropped, but the formula is more intricated. Since we already placed this work in the general setting of stratified spaces, we prefered to present this new application of Morse theory in random topology in this simpler situation. 
		\item In fact, for $N_{\mathbb S^{n-1}}(Z_u)$, we can improve the topological precision: we can impose that the spheres belong to different balls of $M$, in particular we can assume that there cannot be linked. Indeed, there are boundaries of the distinct balls computed for $N_{\mathbb B^n}(E_u)$. 
		\item Note that for $u=0$, all the possible affine topologies have uniform positive densities in the compact algebraic~\cite{gayet2014lower} and Riemannian settings~\cite{GWuniversal} (see also~\cite{sarnak2019topologies} and~\cite{canzani2019topology}), see paragraph~\ref{history}.
	\end{enumerate}
\end{remark}

\paragraph{Betti numbers.}
Morse theory allows us to obtain estimates for the other Betti number $b_i$, 
where for any subset $A\subset \R^n$, 
$ b_i(A) = \dim H_i(A,\R)$
and $b(A) =\sum_{i=0}^n b_i(A)$.
\begin{theorem}\label{corobetti}
	Under the hypotheses of Theorem~\ref{sph0}, assume that $M$ satisfies the further condition~(\ref{bettistrat}) (mild local homology). Then, there exists $c>0$ such that 
	$$\forall u\in \R, \ 
	\mathbb E b(E_u(M,f))= \mathbb E b_0(E_u) \left(1+O_{u\to +\infty} (e^{-cu^2})\right). $$
\end{theorem}
\begin{remark} Theorem~\ref{corobetti} should be true for the nodal set $Z_u$ instead of the excursion set $E_u$, but the proof would involve tedious algebraic topological complications. 
\end{remark}

\paragraph{A refinement.}
If we add a further condition for $M$, namely to be a \emph{locally convex cone space}, see Definitions~\ref{conespacedef} and~\ref{locconv}, and if we use the main result of~\cite{adler}, we can improve Theorem~\ref{sph0} in two ways: a more precise asymptotic and a better bound, but only for $N$, with the notable exception of the class of closed manifolds, see Corollary~\ref{cororef}.
	\begin{theorem}\label{sph2z}Let $\widetilde M$ be a $C^3$ manifold of dimension $n\geq 1$,  $M\subset \widetilde M$ be a compact locally convex $C^2$  cone space of dimension $n$ satisfying condition~(\ref{grass}) (very gentle boundaries), 
		$ \widetilde f: \widetilde M \to \R$ be a random centered Gaussian field satisfying conditions (\ref{C1}) (regularity), (\ref{C2}) (non-degeneracity) and (\ref{C4}) (constant variance), $f=\widetilde f_{|M}$ and $g$ be the metric induced by $\tilde f$ and defined by~(\ref{metric}). 	
		Then, there exists $c>0$ such that 
		\begin{equation}\label{cricri}
		\forall u\in \R, \ \mathbb E N(E_u(M,f)) =
		\sum^n_{k=0}\frac{1}{\sqrt{2\pi}^{k+1}}  \mathcal L_k H_{k-1} (u)e^{-\frac{u^2}2}\left(1+O_{u\to+\infty}(e^{-cu^2})\right),
		\end{equation}
		where the constants $(\mathcal L_k)_k$ are defined below by~(\ref{killing}) and 		$(H_{k})_k$ denote the Hermite polynomials, see~(\ref{hermite}). The error term, included $c$, depends only on the 4-jet of the covariance on the diagonal $\widetilde M\times \widetilde M$.
\end{theorem}
\begin{example} All the examples of Example~\ref{exmo} are locally convex cones, except the last one given by Figure~\ref{tore-bizarre} which is a cone but which is not locally convex.
	\end{example}
\begin{remark} 
	\begin{enumerate}
			\item As for Theorem~\ref{sph0}, a quantitative version of Theorem~\ref{sph2z} is provided by Theorem~\ref{sph2}.
		\item
	It is very likely that the first assertion of Theorem~\ref{sph2z} is true for $N_{\mathbb B^n},$ see Remark~\ref{tark}.
	\end{enumerate}
	\end{remark}
\begin{corollary}\label{cororef}Let $n\geq 1$ be an integer and $M$ be a compact $C^2$ manifold of dimension $n$ with or without boundary, and $f: M\to \R$ be a random centered Gaussian field satisfying conditions (\ref{C1}), (\ref{C2}) and (\ref{C4}). Then, there exists $c>0$ such that (\ref{cricri}) holds. 

If $M$ is a closed manifold, then (\ref{cricri}) writes
\begin{equation}\nonumber
\forall u\in \R, \ \mathbb E N(E_u(M,f)) =
\frac{1}{\sqrt{2\pi}^{n+1}}  \vol_g(M) H_{n-1} (u)e^{-\frac{u^2}2}\left(1+O_{u\to+\infty}(e^{-cu^2})\right).
\end{equation} 

Moreover, again if $M$ is closed, this estimate also
 holds  for 
$N_{\mathbb B^n}(E_u)$, 	$N_{\mathbb S^{n-1}}(Z_u)$ and $N(Z_u)$ instead of 
$ N(E_u).$  
\end{corollary}

\paragraph{Nazarov-Sodin constant.}
 For affine stationnary fields, see condition~(\ref{C3}), the quantitative version of Theorem~\ref{sph0} implies the following simple asymptotic for the constant $c_Z(u)$ defined by Nazarov and Sodin in~\cite{nazarov2}:
	\begin{equation}\label{coco}
	c_{Z}(u)\underset{u\to +\infty}{\sim}
	\frac{1}{\sqrt{2\pi}^{n+1}}\sqrt{\det d^2_{x,y} e(0)}  u^{n-1}e^{-\frac{1}2u^2}.
	\end{equation}
	Here $e$ denotes the covariance function of the field $f$, that is $$ \forall (x,y)\in (\R^n)^2,\ e(x,y)= \mathbb E (f(x)f(y)).$$
	Roughly speaking, $c_Z(u)$ is the volume density of the number of connected components of $Z_u(\R^n, f)$. 
In fact, using the quantitative refinement of Theorem~\ref{sph0} given by Theorem~\ref{sph2}, we obtain a more precise asymptotic with a better error term:
\begin{theorem}\label{ns2}
	Let $f : \R^n\to \R$ be a centered Gaussian field satisfying conditions (\ref{C1})(regularity), (\ref{C2})(non-degeneraticity), (\ref{C3})(stationarity) and~(\ref{C5}) (ergodicity),  and $c_Z(u)$ be the constant defined by Theorem~\ref{nst}. Then, there exists $c>0$ such that
	$$c_{Z}(u)=  \frac{1}{\sqrt{2\pi}^{n+1}}
	\sqrt{\det d^2_{x,y} e(0)}H_{n-1} (u)e^{-\frac{1}2u^2}(1+O_{u\to +\infty}( e^{-cu^2})) ,
	$$
where $H_{n-1}$ denotes the $(n-1)$-th Hermite polynomial given by~(\ref{hermite}) and the constant involved in the error term  
depends only on the 4-jet of $e$ at $0$.
\end{theorem}
\begin{remark}
	\begin{enumerate}
		\item For $n=2$, the asymptotic (\ref{coco}) is a consequence of Swerling's estimate~(\ref{swerling}) given below.
\item Note that Theorem~\ref{ns2} provides the first asymptotic estimate for these enigmatic constants in higher dimensions.  
\item The important case $u=0$ remains unkwown.
\item Note that the equivalent constant $c_E(u)$ for $E_u(\R^n, f)$, instead of $Z_u$, has been proven to exist for $n=2$, see \cite[Theorem 1.3]{beliaev}. Theorem~\ref{ns2} is true for this constant $c_E$ instead of $c_Z$. Notice that it is very likely that $c_E(u)$ is well defined for higher dimensions. 
\end{enumerate}
\end{remark}

\begin{example}\label{bf} Recall that 
	$d^2_{x,y} e(0) = (\partial_{x_i,y_j} e(0))_{1\leq i,j\leq n},$ and let $c_n=\sqrt{\det d^2e (0)}.$ 
	\begin{itemize}
\item	Bargmann-Fock: for $e(x,y)=\exp(-\frac{1}2\|x-y\|^2)$, $c_n=1$.
\item \label{rw} Random waves: for $\displaystyle e(x,y)=\displaystyle 
	\frac{J_{\frac{n-2}2}(\|x-y\|)}{\|x-y\|^{\frac{n-2}2}}$,
	$c_n=
 n^{-\frac{n}2}.
	$
\item \label{frw} Full spectral band random waves: for $ \displaystyle  e(x,y)=
	\frac{J_{\frac{n}2}(\|x-y\|)}{\|x-y\|^{\frac{n}2}}$,
	$c_n=
(n+2)^{-\frac{n+2}2}.
	$
	\end{itemize}
\end{example}

\paragraph{Spin glasses.}
Finally, constructive Morse theory can be applied the so-called~\emph{p-spin spherical spin glass model}, in the context of~\cite{auffinger2013random}.
In this case for any integer $n\geq 1$, $M= \sqrt n \mathbb S^{n-1}\subset \R^n$ and for any integer $p\geq 2$,
the Gaussian random field is defined by
\begin{equation}\label{spinglass}
\forall x=(x_1, \cdots, x_n)\in \sqrt n \mathbb S^{n-1}, \  f_{n} (x)= \sum_{i_1, \cdots, i_p=1}^n a_{i_1\cdots i_p} x_{i_1}\cdots x_{i_p},
\end{equation}
where  the coefficients $(a_{i_1\cdots i_p})_{i_1,\cdots, i_p}$ are independent centered standard Gaussian random variables.
The covariance $e$ of $f_{n}$ satisfies
$$ \forall x, y \in \sn, \,
e(x,y) = n^{1-p}  \langle x,y\rangle^p.$$
Here, the regime consists into increasing the dimension $n$, and looking at the asymptotic behaviour of the sojourn sets $S_{nu}$ (symetric asymptotics hold for $E_{nu}$, see Remark~\ref{sueu}).
\begin{theorem}\label{spin}For any integer $n\geq 2$, 
	let $f_n: \sqrt n \mathbb S^{n-1}\to \R$ be the Gaussian field defined by~(\ref{spinglass}). Then, 
		$$ \forall u<- 2\sqrt{\frac{p-1}{p}},\ 
		\lim_{n\to +\infty}	\frac{1}n \log \mathbb E N_{\mathbb B^{n-1}}\left(S_{nu}(\sqrt n \mathbb S^{n-1}, f_n)\right)= \Theta_{0,p}(u),$$
		where $\Theta_{0,p}$ is the function defined 
		by~\cite[(2.16)]{auffinger2013random}.
The same holds for $N(S_{nu})$ and $N_{ \mathbb S^{n-2}}(Z_{nu})$ instead of  $N_{\mathbb B^{n-1}}(S_{nu})$.
\end{theorem}

\paragraph{The assumptions on the field.}
We now describe the natural assumptions for $\tilde f$  needed in Theorems~\ref{sph0} and~\ref{sph1}. Let $M$ be a Whitney stratified manifold in a manifold $\widetilde M$, see Definition~\ref{stratified}, with local coordinates over each stratum $(x_i)_{1\leq i\leq j}$, and $f : \widetilde M\to \R$ be a centered Gaussian field. The reader only interested in the case where $M$ is a manifold can assume $\widetilde M=M$. 
\begin{enumerate}[series=condition]
	\item \label{C1} (Regularity) The covariance $e : \widetilde M\times \widetilde M\to \R$ is $C^8$ in the neighborhood of $M^2$.
	\item \label{C2} (Non-degeneracity) 
	For any $j\in \{0, \cdots, N\}$, any $x\in \widetilde M$ and any coordinates $(x_i)_{i\in \{1, \cdots, n\}}$, 
	the joint distribution of $$(\partial_i f(x), \partial^2_{kj} f(x))_{i, k\in \{1, \cdots n\}}$$ is non-degenerate.
	\item \label{C4} (Constant variance) The variance of $f$ is constant equal to one, that is $\forall x\in \widetilde M, \ e(x,x)=1.$ 
	\item \label{C3} (Stationarity) If $\widetilde M= \R^n$ for $n\geq 1$, the covariance $e$ is invariant under translations, that is
	$$ \forall (x,y)\in (\R^n)^2, \ e(x,y) = e(x-y,0).$$
	\item \label{C5} (Ergodicity) Under the hypotheses of Condition~(\ref{C3}),
	$e(x,0) \underset{\|x\|\to \infty}{\to} 0.$	
\end{enumerate}
\begin{remark}By Kolmogorov's theorem in~\cite{nazarov2}, 
	Condition (\ref{C1}) implies that the field is almost surely $C^3$, so that the weaker condition of~\cite[(11.3.1)]{adler} is satisfied in coordinates.  
As said in Remark~\ref{rem}, condition~(\ref{C4}) could be dropped, but the formulas are more involved. Moreover, this is a consequence of Condition~(\ref{C3}).
Condition~(\ref{C5}) is only used in Theorem~\ref{ns2} and implies that the action of translations is ergodic.  
\end{remark}
The assumptions for the stratified set need more definitions and results, hence will be defined later in Section~\ref{proof}.

\subsection{Related results}\label{history}

\paragraph{Connected components and critical points.} It seems that the first study of statistics of the number of connected components $N(E_u(M,f))$ of excursion set or $N(Z_u)$ of a random function $f$ in dimension two is due to P. Swerling~\cite{swerling}, in a context of geomorphology. In particular, the author gave lower and upper bounds for the mean number of connected components~\cite[equation (36)]{swerling} of these excursion sets, using Morse-like ideas and estimates of the number of random critical points of given index (maxima, minima and saddle points). The latter study of critical points of random functions in dimensions larger or equal to one began at least in the paper of M. S. Longuet-Higgins~\cite[equation (58)]{longuet}, in a context of oceanography. 

\paragraph{Origins of Morse theory} The idea of linking critical points and topology, which is called now \emph{Morse theory}, can be drawn back to the beautiful and forgotten 1858 article~\cite{reech} by the physicist Fr\'ed\'eric Reech, who computed there the first Morse Euler characteristic using the topology of level lines of the altitude on the Earth, a theorem which A. F. M\"obius generalized~\cite{mobius1863theorie} (citing Reech). Then J. C. Maxwell reproved in 1870 in~\cite{maxwell1870hills}, seemingly unaware of Reech's and M\"obius works.

\paragraph{Euler characteristic.}  
In  1976, the Euler characteristics of the random excursion sets began to be studied~\cite{Adler1976} by R. J. Adler and A. M. Hasofer. 
Note that this invariant is directly accessible \emph{via} Morse theory and critical points, or by Gauss-Bonnet-type formulas, which are local, so that closed formulas can be established through Kac-Rice formulas, on the contrary to the number of connected components.  For spin glasses, more precisely for isotropic Gaussian random fields over $\mathbb S^n$, the study of the Euler characteristics has been done when the dimension $n$ goes to infinity~\cite{auffinger2013random}. Although we won't use it in this paper, it is worth mentionning~\cite{estrade16}, where a central limit theorem was proven for the Euler characteristic of $\chi(Z_0(M,f))$ over larger and larger affine cubes $M$.
On real algebraic manifolds and for random real polynomials, S. S. Podkorytov on the sphere and then T. Letendre in a general setting~\cite{letendre} gave the asymptotic of $\mathbb E \chi(Z_0)$.  
For the proof of the most precise theorem of this article, see Theorem~\ref{sph1}, we use a general asymptotic by R. J. Adler and J.E. Taylor  of $\mathbb \chi(E_u)$ for cone spaces, see Theorem~\ref{euler}.

\paragraph{Large deviations for $N(Z_0)$.} 
In 2006, a regain of interest in connected components was triggered by the work~\cite{nazarov2009number} by F. Nazarov and M. Sodin, who proved that in the context of random eigenfunctions $f$ of the Laplacian over the round $2$-sphere $\mathbb S^2$, $N(Z_0(\mathbb S^2,f))$ has a precise statistics for large eigenvalues. In particular, the average number of $N(Z_0)$ is asymptotic to $cL$, where $c>0$ and $L$ is the increasing eigenvalue. They also proved a large deviation phenomenon. In 2011, the authors of~\cite{gayet2011exponential} proved that for $M$ being a real algebraic surface and $P$ being a random polynomial of large degree $d$, the probability that $N(Z_0(X,P))$ is maximal decreases exponentially fast with $d$ (see also~\cite{diatta2018low} and~\cite{ancona} for recent generalizations and~\cite{rivera2019quasi} for affine fields ; see also~\cite{beliaev2019fluctuations} and~\cite{nazarov2020fluctuations} for estimates of the variance of $N(Z_0)$).
 This work was influenced by former works in random complex algebraic geometry~\cite{shiffman1999distribution}. Note that the latter and \cite{nazarov2009number} were inspired by quantum ergodicity and Berry's conjecture.

\paragraph{Betti numbers of $Z_0$.} 
Non-explicit (like $o(d^n)$) upper bounds and then explicit ones (like $c_n\sqrt d^n$) for $\mathbb E b_i\left(Z_0(X,P)\right)$ were given in the algebraic context in~\cite{gwcrelle} and~\cite{GWBetti}. The authors used Lefschetz and then Morse theory, counting "flip points" of given index, where the random zero set is tangent to a given fixed distribution of hyperplanes.  As said before, this trick that was already used (unknown to the authors) in~\cite{swerling} in dimension 2 for the number of components (then the flip points have index zero or one). 
When large dimension $n$ are studied, large deviations happen for the mean number or critical points of various indexes: the proportion of critical points of indexes close to $n/2$ tends exponentially fast to one~\cite[Theorem 1.6]{GWBetti}. Since the Morse-Euler characteristic equals the Euler characteristic of $Z_0$, weak Morse inequalities (see Thereom~\ref{milnor} assertion~\ref{weak}.) indicate that 
the middle Betti numbers are preponderant compared to the other ones, a phenomenon which is already visible numerically in dimension $n=3$, see~\cite[Figure 5.]{pranav2019}. 
In~\cite{Lerario} A. Lerario and E. Lundberg proved that on the sphere, the mean of $N(Z_0(\mathbb S^n, P))$ has a lower bound growing like $\sqrt d^n$, in various symmetric models. 
In~\cite{gayet2014lower} and~\cite{GWuniversal},  explicit lower bounds for the Betti numbers were given in algebraic and Riemannian settings. In a different spirit, \cite{lerario2016gap} dealt with mean Betti numbers of random quadrics with increasing dimension.

\paragraph{Diffeomorphism type of $Z_0$.} 
In 2014, the diffeomorphism type of the random nodal sets $Z_0$ began to be studied in~\cite{gayet2014lower}. The authors proved that for $X$ being an $n$-dimensional real algebraic manifold and $P$ being a random polynomial of degree $d$, for any 
affine compact hypersurface $\Sigma\subset \R^n$,  the average $\mathbb E N_{\Sigma}(Z_0(X,P))$ of the components of $Z_0$ diffeomorphic to $\Sigma$ also grows at least like $c_\Sigma \sqrt d^n$, where $c_\Sigma>0$ can be made explicit.  The same was then proven in~\cite{GWuniversal} for a random sum of eigenfunctions of Laplacian for eigenvalues up to a large increasing number $L$. 

\paragraph{Asymptotic values.} 
In 2016, F. Nazarov and M. Sodin proved in~\cite{nazarov2} that in a very general context, for stationary  affine Gaussian random  field,  $\frac{1}{\vol(B^n_R)}\mathbb E N(Z_0(B_R^n, f)) $ converges to a positive constant $c_Z(0)$ when $R$ grows to infinity, see Theorem~\ref{nst} below. 
In 2019, P. Sarnak and I. Wigman, and P. Sarnak and Y. Canzani,  gave a version of this result in ~\cite{sarnak2019topologies} and~\cite{canzani2019topology}  for the number $N_{\Sigma}(Z_0)$ defined above, in the Laplacian context. In~\cite{wigman19}, I. Wigman gave a version of the Nazarov-Sodin asymptotic for Betti numbers of the components of $Z_0(B_R^n, f)$ which do not intersect the boundary of $B_R^n.$ Note that in the contrary to $i=0$, it could happen that for $i\geq 1$,  a unique large connected component of $Z_u(B_R^n)$ has a large $b_i$ and touches the boundary.

\paragraph{Estimates in dimension 2.} 

The values of the average of the numbers $N(Z_0)$, $N(E_u)$, $N_{\Sigma}( E_u)$ or their asymptotics $c_Z(0)$ or $c_E(u)$ are unkwnon, and the known bounds for them are related to either critical points, which are far easier to compute, or the barrier method (see~\cite{nazarov2009number}). 
Until the present work, the dimension $n=2$ was the only case where asymptotics has been done. Indeed in the affine case and for  isotropic smooth centered Gaussian fields,  \cite[Equation (36)]{swerling} implies that
\begin{equation}\label{swerling}
\forall u\in \R, \ 
\frac{ \mathbb E N(Z_u(B_R,f))}{\sqrt{\det d^2_{x,y}e (0)}\vol (B_R)}
=
\frac{1}{\sqrt{2\pi}^3}u   e^{-\frac{1}2u^2} 
\left(1+O_{u\to +\infty}(
 e^{-cu^2})
\right)
\end{equation}
where $c>0$ depends only on the 4-jet of $e$ at $0$ (see also~\cite[Corollary 1.12 and Proposition 1.15]{beliaev} for non isotropic fields).
Since $H_1=\text{Id}$, this estimate is the same as our Theorem~\ref{ns2} for $n=2$.
Also in dimension $2$, T. L. Malevich gave bounds for $\mathbb E  N(Z_0)$ in~\cite{malevich72},
for fields with positive correlations. The method is more direct there than in Swerling's paper, but less precise. In~\cite{ingremeau2018}, the authors gave an explicit lower bound for $c_Z(0)$ in the case of planar random waves. 

\paragraph{Estimates in higher dimensions.} 
In higher dimensions, Nicolaescu~\cite[Theorem 1.1]{nicolaescu2015critical} gave a universal upper bound for $\mathbb E N(Z_0)$ in the Riemannian setting, using the number of local minima. As said before explicit lower bounds for $\mathbb E N(Z_0)$ were given by~\cite[Corollary 1.3]{gayet2014lower} and~\cite[Corollary 0.6]{GWuniversal}, and for  $\mathbb E N_{\Sigma}(Z_0)$ for  $\Sigma$ being the sphere or products of spheres (in order to get higher Betti numbers), in compact algebraic and Laplacian (even elliptic operators) contexts. For instance, if $(M,g)$ is a compact Riemannian manifold and $f$ is a random sum of eigenfunctions of the Laplacian with eigenvalue bounded above by $L$ (\cite[Corollary 0.6]{GWuniversal} and \cite[Corollary 0.3]{GWasian}), for $L$ large enough and for any $i\in \{0,\cdots, n-1\}$,
\begin{equation}\label{sand}
	\exp(-e^{312 n^{3/2}})\leq 
	\frac{\mathbb E N_{S^i\times S^{n-1-i}}(Z_0)}{\sqrt L^n\vol_g (X)} \leq
\frac{\int_{\sym(i,n-1-i,\R)}|\det A| d\mu(A)}{\sqrt \pi^{n+1}\sqrt{(n+2)(n+4)^{n-1}}}
\leq  \exp(-c_i n^2)
\end{equation}
where $d\mu(A)$ is an explicit  universal measure (which is for instance GOE in the algebraic version) and $c_i>0$. The upper bound in~(\ref{sand}) given the average of critical points is likely to have the good order since the Morse-Euler characteristic equal the topological one.  Note also that all of these estimates should be essentially true with  similar results in the affine case, at least for Bargmann-Fock and full band random waves, see Example~\ref{bf} below for the definition, since the compact case converges, after rescaling at order $1/\sqrt d^n$ or $1/\sqrt L^{n/2}$ near a fixed point, to the affine model. Note also that in principle, these estimates
for the nodal hypersurfaces $Z_0$  should be  adapted for the topology of the level set $Z_u$ or $E_u$ with non-zero $u$. 

\paragraph{Bumps and Euler characteristic.} 
In~\cite{adler}, it is given a closed formula 
the mean Euler characteristic 
of $E_u$, see Theorem~\ref{euler} below. 
Since the Euler characteristic of a ball is one, under the belief that most of components of $E_u$ are balls, the estimate for the Euler characteristic given by Theorem~\ref{euler}  should be true for the mean number of connected components $\mathbb E  N(E_u)$, and even for the mean number $\mathbb E N_{\mathbb B}(E_u)$ of components diffeomorphic to a ball.
We prove that it is true, see Theorem~\ref{sph1}.  
As a final remark, note that the shape near a non-degenerate local maximum is automatically a  topological ball.  It has been proven in~\cite[(6.2.12)]{adler10} that the shape of this ball is quite precise, but it does not allow to estimate the number of connected components as in our Theorem~\ref{sph0}.

\paragraph{Random topology and cosmology.}
The topology of random excursion or nodal sets, in particular Betti numbers, has become 
a subject of interest in cosmology in the last decade, at least since~\cite{van2011}. Here, the field is the mass density. We refer to the survey~\cite{pranav2019}
for references to this subject. 

\paragraph{Questions} We finish this section with questions that the present work arises. 
\begin{enumerate}
	\item Is there a closed formula for $\mathbb E N(E_u)$ and its various avatars, at least for affine isotropic fields? 
	\item In particular, what is the value of $c_Z(0)$?
	\item Is it possible to obtain an asymptotic of $\mathbb E N_{\Sigma}(E_u)$ for other manifolds $\Sigma$ than $\mathbb B^n$?
	\end{enumerate}

\paragraph{Structure of the article.} 
\begin{itemize}
\item Section~\ref{morsem} is of deterministic nature and recalls the classical main elements of Morse theory on a compact smooth manifold without boundary, that is how the topology of the sublevel (sojourn set) $S_u(M,f)$ of a Morse function $f$ changes when passing a critical point. Since the change of topology of the level set (nodal set) $ Z_u(M,f)$ is in general not treated, we provide the results we need in the sequel. 
In this section we treat the spin glass situation, since the comparison of the average number of critical points of various indices has been already done in~\cite{auffinger2013random}. 
\item 
Section~\ref{randomf} handles with random Gaussian fields  over manifolds and their critical points. We follow the elegant stochastico-Riemannian setting developped in~\cite{adler}, where the metric is induced by the random field. Then we compute the average number of critical points of given indices in the spirit of~\cite{adler10}, where it was proven that, for large $u$, local maxima predominate exponentially fast amongst critical points of $f$ in $E_u$. But here we need and provide a expanded version of it with explicit error terms and for manifolds. 
\item 
Section~\ref{morsestrat} is of deterministic nature and explains the general setting, that is the definition of Whitney stratified sets. It then presents
the main features of \emph{stratified Morse theory},
a vast generalization of classical Morse theory developped in the book~\cite{goresky}. 
\item 
Section~\ref{proof} is devoted to the proofs of the theorems for the general Whitney stratified spaces.  
We then explain how the full result of~\cite{adler} for the mean Euler characteristic can be used for $\mathbb E N(E_u)$ when we assume that the Whitney space is in fact a locally convex cone space. We then prove the asymptotic formula of the Nazarov-Sodin constant $c_Z(u)$. 
\end{itemize}
\paragraph{Acknowledgements.} The author thanks Antonio Auffinger for his valuable expertise about~\cite{auffinger2013random} and the first part of  the proof of Theorem~\ref{spin}. He also thanks Fran\c{c}ois Laudenbach for the part of Remark~\ref{tark} concerning manifolds with boundary. 
The research leading to these results has received funding from the French Agence nationale de la recherche, ANR-15CE40-0007-01 (Microlocal) and ANR-20-CE40-0017 (Adyct).

\section{Classical Morse theory}\label{morsem}

\subsection{Change of sojourn sets}\label{morsemm}

Classical Morse theory~\cite{milnor} is a way to understand part of the topology of a compact smooth manifolds using one $C^2$ function on it, as long as its critical points are non-degenerate, that is its Hessian at these points are definite. As said in the introduction, it seems that the first occurrence of this circle of ideas draws back to F. Reech~\cite{reech}. Let $M$ be a compact smooth $n$-manifold without boundary,  $f : M\to \R$ be a $C^2$ map and let $u\in \R$. Recall the definition of the sojourn set 
 $S_u(M,f) = \{x\in \R, f(x)\leq u\}.$
\begin{remark}\label{sueu}
	\begin{enumerate}
\item In the Morse theory tradition,  $S_u(M,f)$ is written $M_u$ or $M_{\leq u}.$ Since $f$ will be random and hence will change, we prefer the probabilistic notation.
\item Note that 
$ S_u(M,f) = E_{-u}(M,-f).$ Hence, for centered fields, the law of the subsets $S_u(M,f)$ is the same of the one of the subsets $E_{-u}(M,f)$, so that in particular 
$$ \forall u\in \R, \ \mathbb E S_u(M,f) = \mathbb E E_{-u}(M,f) .$$
\end{enumerate}
\end{remark}
Recall that a \emph{critical point} $p$ of $f$ is a point of $M$ such that 
$df(x)=0$. At a critical point, we can define the second differential $d^2 f(x)$ in any coordinate system. Then, $d^2f(x)$  has a definite signature independent of the coordinates. Define
\begin{align}\label{index} 
\forall  A\in \sym(n,\R), \ \Ind(A)&=\# \text{Spec}(A)\cap \R^-,
\end{align}
where $\sym(n,\R)$ denotes the set of real symmetric matrices of size $n$ and $\text{Spec}$ the spectrum.
For any subset $E\subset M$ and $i\in \{0,\cdots, n\}$ and $f:M\to \R$ any  Morse function, let
\begin{eqnarray*}
	\crit_i(E,f)=\{x\in E, df(x)=0 \text{ and } \Ind(d^2 f(x))=i\},\\
	C_i(E,f)=\# \crit_i(E,f) \text{ and }
	C(E,f)= \sum_{i=0}^{\dim M} C_i(E,f).
\end{eqnarray*}
We will omit $f$ when it is tacit. 
Let for any integer $i$, 
 $$b_i(S) = \dim H_i (S, \R)$$ be the $i$-th Betti number of $S$. In particular, $b_0(S)$ is the number of connected components of $S$.
\begin{definition}\label{morsedef} Let $M$ be a $C^2$ manifold and $f : M\to \R$  be a $C^2$ function.
The map $f$ is said to be \emph{Morse} if 
the critical points of $f$ 
are isolated, and non-degenerate. The latter means that is its Hessian in coordinates is definite. 
\end{definition}
This definition does not depend on the chosen coordinates.

\begin{definition} 
Let $B\subset A$ and $S_-\subset S_+ $ be topological spaces and $g: B \to S_-$ be a continuous map,  such that the identity
	map $S_- \subset S_-$ extends to a homeomorphism
	$$S_+\sim S_-\cup_g A=S_-\sqcup A/\sim ,$$
	where for all $y\in S_-$ and $x\in B, y\sim x$ whenever $y=g(x)$.
	Then we shall say that $S_+$ \emph{is obtained from $S_-$ by attaching the pair} $(A, B)$
	and we will write
	$$S_+ = S_-\cup_g (A, B).$$
	\end{definition}
Note that when $B=\emptyset$, then $S_+=S_-\sqcup A$. For $n\in \Nn$, we will use the notation $D^n$ for the unit ball $\mathbb B^n \subset \R^n$; to be clear, $D^0$ is a point. By $\partial D^n$ we denote the sphere $\mathbb S^{n-1}\subset \R^n$, so that $\partial D^1 \{-1, 1\}$ and $\partial D^0=\emptyset$.
We now sum up the main features in classical Morse theory we will use:
\begin{theorem}\label{milnor}
	Let $M$ be a compact smooth manifold of dimension $n$ and $f : M\to \R$ be a Morse function. Then, the following holds:
\begin{enumerate}
	\item (Invariance between two non-critical values)\label{homot}~\cite[Theorem 3.1]{milnor} Let $u\leq v\in \R$ be such that $[u,v]$ does not contain any critical value of $f$. 
	Then $S_v(M,f)$ is diffeomorphic (up to boundary) to $S_u(M,f)$.
\item (Change at a critical point)\label{matsu} \cite[Proposition 4.5]{goresky}  Let $u\in \R$ be such that there is a unique critical point $p\in M$ in $Z_u(M,f)$, and assume that $p$ has index $i\in \{0,\cdots, n\}$. Then, for any $\epsilon>0$ small enough, 
		the manifold with boundary  
		$S_{u+\epsilon}(M,f)$ is homeomorphic to $$S_{u-\epsilon}(M,f)\cup_g (D^i\times D^{n-i},\partial D^i\times D^{n-i}),$$ 
		where $g  : \partial D^i\times D^{n-i}\to Z_{u-\epsilon}$  the attaching map is an embedding. In particular their boundaries are homeomorphic.
			\item (Components diffeomorphic to a ball)\label{nball} Let $u\in \R $ be a non-critical value of $f$. Then, any connected component of $S_u(M,f)$ containing a unique local minimum and no other critical point is diffeomorphic to a $n$-ball $\mathbb B^n$. 
\item\label{krilling} (Killing of a component) Under the hypotheses of assertion~\ref{matsu}., for any small enough positive $\epsilon$,
$$b_0(S_{u+\epsilon})\geq b_0(S_{u-\epsilon}) - {\bf 1}_{\{p\in \crit_1(M,f)\} }.$$
	\item\label{weak} (Weak Morse inequality)~\cite[Theorem 5.2]{milnor} For any non-critical level $u\in \R$, $$\forall i\in\{0,\cdots, n\},\,  b_i(S_u (M,f))\leq C_i(S_u).$$
	\item~\label{eulermorse} (Morse Euler characteristic)~\cite[Theorem 5.2]{milnor} For any non-critical level $u\in \R$,
	\begin{equation}\label{chi}
	 \chi(S_u) = \sum_{i=0}^n (-1)^i C_i (S_u).
	 \end{equation}
	\end{enumerate}
\end{theorem}
As explained in the introduction, F. Reech proved~(\ref{chi}) in~\cite{reech} for $M=\mathbb S^2$.

	\begin{figure}
	\centering
		\includegraphics[width=0.47\textwidth]{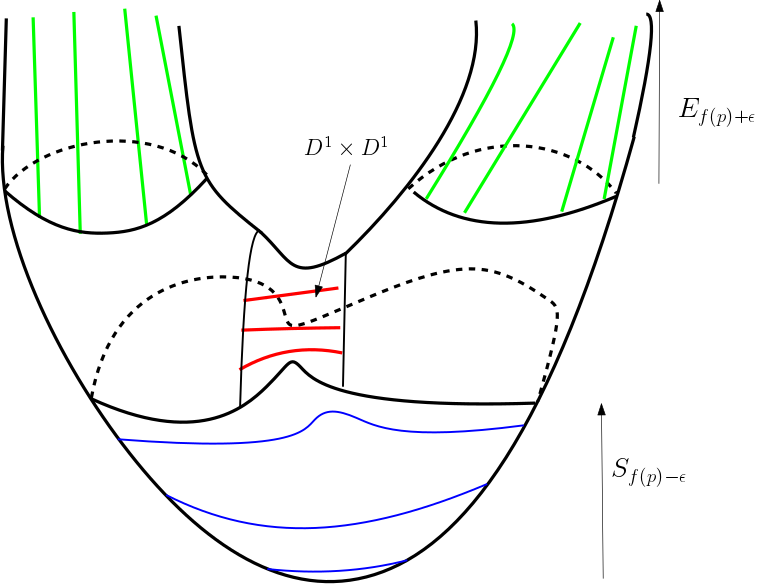}
	\includegraphics[width=0.47\textwidth]{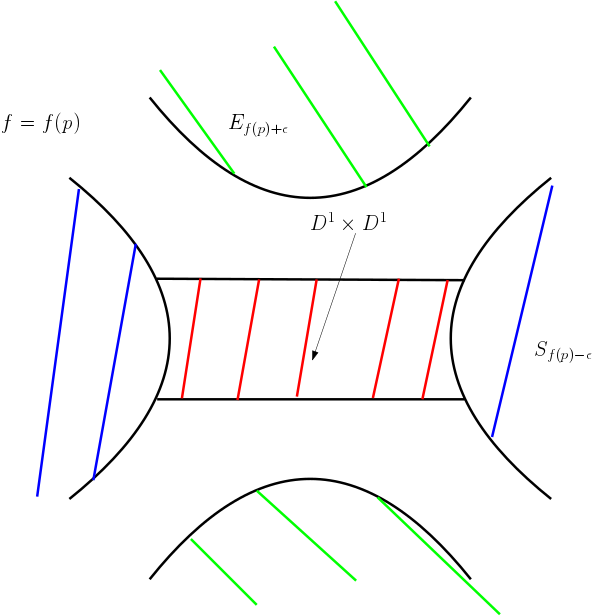}
	\caption{Attaching a handle to $S_{u-\epsilon}$. In this drawing, $n=2$ and the index equals $i=1$. The sojourn set $S_{u+\epsilon}$ (the complementary of the green shaded part of the surface) is homeomorphic to $S_{u-\epsilon}\cup_g (D^1\times D^1, \partial D^1\times D^1)$.  }
	\label{anse}
\end{figure}
\begin{corollary}\label{coromilnor} Under the hypotheses of Theorem~\ref{milnor},
	for any non critical real $u\in \R$,
	\begin{equation}\label{inek}
	 0\leq C_0(S_u(M,f))-N(S_u)\leq C_0(S_u)-N_{\mathbb B^n}
	(S_u)\leq
	 \sum_{i\geq 1} C_i(S_u).
	 \end{equation}
\end{corollary}
\begin{preuve} The first inequality is a consequence of Theorem~\ref{milnor} assertion~\ref{weak}. for $i=0$. The second one is trivial. For the last one, 
	by assertion~\ref{nball}. any critical point with vanishing index produces  component homeomorphic to a ball, and a critical point of index larger or equal to 1 can change the topology of at most one component of $S_u$, hence the last inequality in~(\ref{inek}). 
	\end{preuve}

\subsection{Change of nodal sets}

The second assertion of Theorem~\ref{sph0} about the nodal sets $Z_u$ needs the following further result which we could not find in the litterature (see however~\cite{knauf} for Betti numbers estimates).
\begin{proposition}\label{morseZ}
	Let $M$ be a compact smooth manifold without boundary and $f : M\to \R$ be a Morse function. Then, 
	\begin{enumerate}
		\item	(Invariance)	For any pair of reals $u<v$ such that $f$ has no critical value in $[u,v]$,  $Z_u(M,f) $ is diffeomorphic to $Z_v(M,f)$. 
	\item\label{taka} (Change at a critical point) For any $u\in \R$, if $p$ is a unique  critical point in $Z_u(M,f)$ with value $u=f(p)$, then 
	for $\epsilon$ positive and small enough,
	$$ |b_0(Z_{u+\epsilon}(M,f)) -b_0(Z_{u-\epsilon})|\leq 6.$$
	\item\label{taka0} (Components diffeomorphic to a sphere) 
	Under the hypotheses of assertion~\ref{taka}, if $p$ has vanishing index, then 
	$$ Z_{u+\epsilon}\sim_{diff} Z_{u-\epsilon} \sqcup \mathbb S^{n-1}.$$
		\end{enumerate}
	\end{proposition}
\begin{preuve}
	The first assertion is a direct consequence of Theorem~\ref{milnor} assertion~\ref{homot}.	The third one is a consequence of Theorem~\ref{milnor} assertion~\ref{nball}. For the second point, assume that the index of $p$ is $i$.
	By Theorem~\ref{milnor} assertion~\ref{matsu},
 $$ S_{u+\epsilon}\sim_{homeo} S_+=S_{u-\epsilon}\cup_g (A,B),$$
where 
	 $$ (A,B)=(D^i \times D^{n-i},\partial D^i \times D^{n-i})
\text{	 and }
	g : B\to S_{u-\epsilon}$$ is the attaching map, which is an embedding. 	In particular,
since $M$ has no boundary, $\partial S_+$ is homeomorphic to $Z_{u+\epsilon}$.
Let
 $$Z_-:= Z_{u-\epsilon}\setminus g(B).$$
 Then $\partial S_+\sim Z_- \cup_{\partial g} (\partial A\setminus B)$, where 
	$ \partial g = g_{| \partial (\partial A\setminus B)}.$
	Put a metric on the handle $A$ and for  $\eta>0$, let 
	$$U=\{z\in \partial A, \text{dist}(z,B)\geq 2\eta\}\text{ and } 
	V = Z_-\cup \{z\in \partial A\setminus B, \text{dist}(z,B)\leq 3\eta\}.$$ For $\eta$ small enough, 
	$$U\sim \partial A\setminus B=D^i \times \partial D^{n-i}\text{ and } V\sim_{retract} Z_-.$$ 
	By Mayer-Vietoris, since $U\cup V = \partial S_+$ and $U\cap V $ is homeomorphic to 
	$ \partial B\sim \partial D^i \times \partial D^{n-i},$
	there exists a long exact sequence
	$$ 
\cdots\to	H_0(\partial B) 
	\overset{\alpha}{\to} H_0(\partial A\setminus B)\oplus H_0(V)\to
	H_0(\partial S_+)\to 0,$$
	so that 
	\begin{equation}\label{amor}
	b_0(\partial S_+)= 
	b_0(\partial A\setminus B)+b_0(V)-\rank (\alpha).
	\end{equation}
	In order to estimate $b_0(V)$, let 
	 $W$ be a small tubular neighborhood of $ g(B)$ in $Z_{u-\epsilon}
	$.
	Note that, since $g(B)\sim_{homeo} B$, 
	$$W\cup V\sim_{retract} Z_{u-\epsilon} \text{ and } W\cap V \sim_{retract} \partial B.$$
Then, again by Mayer-Vietoris,
	$$ 
	H_0(\partial B) 
	\overset{\beta}{\to} H_0(W)\oplus H_0(V)\to
	H_0(Z_{u-\epsilon})\to 0,$$ 
	so that, since $W\sim_{retract} g(B)\sim B$, 
	$$ b_0(V) = b_0(Z_{u-\epsilon})+\rank (\beta) -b_0(B).$$
	Finally, by~(\ref{amor}), we obtain
	\begin{equation}\label{toto}
	b_0(Z_{u+\epsilon})-
	 b_0( Z_{u-\epsilon})= b_0(\partial A\setminus B)+\rank (\beta) -b_0(B)-\rank (\alpha),
	 \end{equation}
	 where $b_0(B)\leq 2$ and
	 $$ \rank (\alpha)\leq b_0(\partial D^i \times\partial D^{n-i})\leq 4\text{ and } \rank(\beta) \leq b_0(\partial D^i \times \partial D^{n-i})\leq 4.$$
	 Hence, $b_0(Z_{u+\epsilon})-
	 b_0( Z_{u-\epsilon})\in [-6,6].$
	 which implies the result. 
\end{preuve}

	The following corollary is analogous to Corollary~\ref{coromilnor}.
\begin{corollary}\label{pasteque} 
	Under the hypotheses of Proposition~\ref{morseZ}, for any non-critical value $u\in \R$, 
	$$|N(Z_u(M,f)) - C_0(S_u)|\leq 6\sum_{i=1}^n C_i(S_u).$$
	The same holds for $N_{\mathbb S^{n-1}}(Z_u)$ instead of $N(Z_u)$.
\end{corollary}
\begin{preuve} This is an immediate consequence of Proposition~\ref{morseZ}.
\end{preuve}

\subsection{Spin glasses}

We finish this section by proving  Theorem~\ref{spin}. In the setting explained in the introduction, two shortcuts happen. First, the field is defined over the sphere, so that we don't need stratified Morse theory. Second, the comparison of the average number of critical points of given index is done in~\cite{auffinger2013random}.  
We begin by recalling the main results of~\cite{auffinger2013random} which we will use. 
\begin{theorem}~\cite[Theorems 2.1, 2.5 and 2.8]{auffinger2013random}
	\label{auf} Let $n\geq 1$ and $p\geq 2$ be integers, and $f_{n}: \sqrt n \mathbb S^{n-1}\to \R$ be	the Gaussian random field defined by~(\ref{spinglass}).	Then, for any index $i\in \{0, \cdots, n-1\}$ and any $u\in \R$,	
	\begin{equation}\label{sissi}
		\mathbb E C_{i}\left(S_{nu} (\sqrt n \mathbb S^{n-1}, f_n)\right)=\sqrt{\frac{8}p}(p-1)^{n/2} \mathbb E_{GOE_n} 
	\\ \left[ e^{-n\frac{p-2}p\lambda_i^2}{\bf 1}_{\lambda_i\leq 
		\sqrt{\frac{p}{2(p-1)}}u}\right],
	\end{equation}
	where the $GOE_n$ measure is the classical measure on the space of real symmetric matrices of size $n$, see~\cite[(2.6)]{auffinger2013random} and $\lambda_0\leq  \cdots\leq \lambda_{n-1}$ denote the (real) eigenvalues of the random symmetric matrix. 
	Moreover, for any $i\in \Nn,$
	$$ \forall u\in \R, \	
	\lim_{n\to +\infty}	\frac{1}n \log \mathbb E C_i\left(S_{nu}(\sqrt n \mathbb S^{n-1}, f_n)\right)= \Theta_{i,p}(u)$$
	where $\Theta_{i,p}$ 
	is defined 	by~\cite[(2.16)]{auffinger2013random}.
\end{theorem}
\begin{preuve}[ of Theorem~\ref{spin}]
Let  $i\in \{1,\cdots n-1\}$ and $u\leq 0$. Then, by~(\ref{sissi}),
	\begin{eqnarray*}	\mathbb E C_i (S_{nu}(\sqrt n \mathbb S^{n-1},f_n)
		) &\leq&  \sqrt{\frac{8}p}(p-1)^{n/2}e^{-n\frac{p-2}{4(p-1)}u^2} \mathbb P_{GOE_n} \left[\lambda_i\leq 
		\sqrt{\frac{p}{2(p-1)}}u\right]\\
		&\leq &	\sqrt{\frac{8}p}(p-1)^{n/2}e^{-n\frac{p-2}{4(p-1)}u^2} \mathbb P_{GOE_n} \left[\lambda_1\leq 
		\sqrt{\frac{p}{2(p-1)}}u\right].
	\end{eqnarray*}
	Let $u\leq -E_\infty<0$, where $E_\infty=2\sqrt{\frac{p-1}{p}}.$
	By the large deviation result given by~\cite[Theorem A.9]{auffinger2013random}, and paying attention to the used normalizations~\cite[Remark 2.4]{auffinger2013random}, $I_1$ is the rate function for $u\leq -E_\infty$. In particular,
	$$\forall u\leq -E_\infty,\  \limsup_{n\to +\infty} \frac{1}n\log  \mathbb P_{GOE_n}
	\left[\lambda_1\leq 
	\sqrt{\frac{p}{2(p-1)}}u\right]= -2I_1(u),
	$$ 
	where $I(1)\geq 0$ is defined by~\cite[(2.13)]{auffinger2013random} and vanishes only at $-E_\infty$.
	The former inequality, the latter limit and the definition of $\Theta_1$ imply that 
	$$ \limsup_{n\to +\infty} \max_{i=\{1, \cdots,n-1\}}\frac{1}n\log 	\mathbb E C_i (S_{nu})\leq 
	\Theta_{1,p}(u).$$
	However,  
	$\forall u\leq -E_\infty$, $\Theta_{1,p}(u)=\Theta_{0,p}(u) -I(1),$ so that
	by Theorem~\ref{auf}, 
	$$\limsup_{n\to +\infty}\max_{i=\{1, \cdots,n-1\}} \frac{1}n\log 	\mathbb E C_i (S_{nu})\leq \lim_{n\to +\infty} \frac{1}n\log 	\mathbb E C_0 (S_{nu})-I_1(u).$$
	Now, fix $u<-E_\infty$. Then, there exists $N\in \Nn^*$, 
	such that 
	$$ \forall n\geq N,\ \mathbb E C_0 (S_{nu})-
	\sum_{i=1}^{n-1}\mathbb E C_i (S_{nu})\geq \mathbb E C_0 (S_{nu})(1-ne^{-nI_1(u)/2}).
	$$
	By Corollary~\ref{coromilnor} and Theorem~\ref{auf},
	this implies that 
	$$\forall u<-E_\infty,\ \liminf_{n\to +\infty}\frac{1}n\log 	\mathbb E N_{\mathbb B^{n-1}} (S_{nu})\geq \Theta_{0,p}(u).$$
	By the weak Morse inequality (Theorem~\ref{milnor}), the latter $\liminf$ is bounded above by $\Theta_{0,p}(u)$ as well, hence the result. 
	The same holds for $N_{\mathbb S^{N-2}}(Z_u)$ instead of $N_{\mathbb B^{n-1}}(E_u)$ applying Corollary~\ref{pasteque} instead of 
	Corollary~\ref{coromilnor}.		
\end{preuve}

\section{Gaussian fields over manifolds}\label{randomf}

\subsection{The induced Riemannian geometry}

\paragraph{Riemannian generalities.}\label{indu}
Let $(M,g)$ be a $C^2$ Riemannian manifold. 
The curvature operator induced by $g$ is denoted here by $R$, which is a two-form over $M$ with values in $TM$~\cite[(7.5.1)]{adler}. It induces the curvature also written $R$~\cite[(7.5.2)]{adler}: 
$$ \forall X,Y,Z,W \in TM,\  
R(X,Y,Z,W) = g(R(X,Y)Z,W).$$
Assume now that $N \subset M$ is a submanifold of codimension at least 1 (later, $N$ will be a stratum $\partial_j M$ of $M$ for $j\geq 1$, or $N$ will be $M$ inside $\widetilde M$). 
The second fundamental form associated to $N, M$ and $g$ is defined by~\cite[(7.5.8)]{adler}:
\begin{equation}\label{second}
\forall X,Y\in TN,\  S(X,Y) = P_{T^\perp N} ( \nabla_X Y),
\end{equation}
where $\nabla $ denotes the Levi-Civita connection associated to $g$, see~\cite[(p.163)]{adler}, and $$P_{T^\perp N} : TM_{|N} \to T^\perp N$$ denotes the orthogonal projection of $TM$ onto the normal bundle $T^\perp N $ over $N$.
Finally, 
for $\nu \in  T^\perp N$, we define
\begin{eqnarray}\label{snu}
S_\nu : TN^2 &\to &\R,\nonumber \\
 \forall X,Y\in TN, \ S_\nu(X,Y) &=& g (S(X,Y),\nu).
 \end{eqnarray}
We also define the first fundamental form:
$$ I : TN\times TN\to \R$$
to be the scalar product $g$, that is $I(X,Y)=g_{|N}(X,Y)=g(X,Y)$. 
Recall also that for $f: M\to \R$ a $C^2$ map, 
$$ \nabla^2 f : TM\times TM \to \R$$
is defined by~\cite[12.2.7]{adler}
$$ \forall X,Y\in \Gamma(TM), \ \nabla^2f(X,Y)= XYf - \nabla_X Yf,$$
where $\Gamma(TM)$ denotes the set of $C^2$ vector fields. Note that at a critical point, $\nabla^2 f (X,Y)=XYf.$ Since $\nabla$ is torsion-free, this is a symmetric bilinear form which depends only on the value of the fields at the point where it is computed.	Recall that for a function $f$, $df=\nabla f$.

\paragraph{Metric induced by a Gaussian field.}

Let $M$ be a $C^2$ manifold and $f: M\to \R$ be a centered Gaussian field
satisfying 
the hypotheses (\ref{C1}) (regularity) and (\ref{C2}) (non-degeneracity for $df$). Recall  that  
$f$ induces a metric over $M$ 
by~\cite[(12.2.1)]{adler}:
\begin{equation}\label{metric}
\forall x\in M, X,Y\in T_x M, \ 
g(X,Y) =\mathbb E(df(x)X, df(x)Y),
\end{equation}
where $df(x)$ denotes the differential of $f$ at $x$.

\begin{example}\label{remstat}
	If $ f: \R^n \to \R$ satisfies condition~(\ref{C3}) (stationarity), that is, there exists a function $ k : \R^n \to \R$, such that is its covariance function $e $ satisfies $\forall (x,y)\in \R^n, \ e(x,y) = k(x-y)$, 
	then the associated metric $g$ satisfies
	$ g = -d^2 k(0).$
	In this case, $\nabla = d$, $R=0$ and $S_\nu = d\nu$.
\end{example}

We present now very useful and elegant computations from~\cite{adler}.
\begin{proposition}\label{gorille}\cite[(12.2.13),(12.2.15)]{adler}
	Let $ \widetilde M$ be a $C^3$ Riemannian manifold, $M\subset \widetilde M$ be a $C^2$ submanifold, $\widetilde f: M\to \R$ be a Gaussian  centered field satisfying conditions~(\ref{C1}) (regularity), (\ref{C2}) (non-degeneraticity) and (\ref{C4}) (constant variance), and $g$ be the metric~(\ref{metric}) associated to $\widetilde f$.
	Then, for any $x\in M$, $\nu^*\in T^*_xM$, $u\in \R$, 
	\begin{eqnarray}
	E:=\mathbb E \left(	\nabla^2 f \ |\  d f = \nu^*, f=u\right)&=&-u I\\
	\forall X,Y\in TM,\  \mathbb E \left(XYf  \ |\  d f=  \nu^*, f=u\right)&=&-ug(X,Y) +\nu^*(\nabla_X Y)\\
	\mathbb E \left(	(\nabla^2 f - E)^2 \ |\  d f = \nu^*, f=u\right)&=&-(2R+I^2),\label{pepette}
	\end{eqnarray}
	where everything is computed at $x$.
	Moreover, the right-hand side of equation~(\ref{pepette}) is the covariance of the conditionned second derivative~\cite[Lemma 12.3.1]{adler}, and the operators are restricted to $T_xM$.
\end{proposition}
In equation~(\ref{pepette}), the square of a 2-tensor $P$ is defined by~\cite[(7.2.5)]{adler}:
$$ \forall X,Y,Z,W\in T_x M,
\ P^2  (X,Y,Z,W)=2\left(P(X,Z)P(Y,W)-P(X,W)P(Y,Z)\right).$$

\subsection{Critical points}\label{criticalpoints}
Critical points are the key elements of Morse theory. Luckily, they are in principle pretty simple to compute in average, because of the Kac-Rice formula. We establish two related results, the idea of the third one coming from~\cite{adler10}.

\paragraph{Setting and notations.}
We begin with the general setting of this part. 
Let $\widetilde M$ be a $C^3$ manifold of dimension $N\geq 1$, $M\subset \widetilde M$ be a $C^2$ submanifold of dimension $n\in \{0,\cdots, N\}$ with compact closure,  $\widetilde f : \widetilde M\to \R$ be a Gaussian field satisfying conditions~(\ref{C1}), (\ref{C2}) and (\ref{C4}), and $f=\widetilde f_{|M}$, and $g$ the metric~(\ref{metric}) induced by $\widetilde f$. 
For any $i\in \{0,\cdots, n\}$, we denote by $C_i(M)$ the number of critical points of $f$ of index $i$.
\begin{itemize}
	\item
Assume first that $n\leq N-1$. 	
For $t\in \R$, $x\in M$ and $\nu\in T_x^\perp M$, define 
$\mu_{t,\nu}$ the Gaussian measure over $\sym(n,\R)$
with average $-tI + S_\nu$ and variance
$ (-2R +I^2)$,
\begin{equation}\label{mutnu}
\forall t\in \R, (x,\nu)\in T_x^\perp M,\ 	\mu_{t,\nu}\sim 
N(-tI + S_\nu, -2R +I^2),
\end{equation}
viewed and restricted to an orthonormal basis for $T_x M$, where $S_\nu$ was defined by~(\ref{snu}). More concretely, for any $x\in M$,
we fix $(E_i)_{i\in T_xM}$ an orthonormal basis of $(T_xM,g)$
so that $(\R^n,g_0)$ is identified with $(T_xM,g)$ through it, as well
as the operators $I$ (which becomes the identity matrix), $R$ and $S_\nu$. Then, if $(Y_{kl})_{1\leq k,l\leq n}\in \sym(n,\R)$ is equipped with 
the measure $\mu_{t,\nu}$, then  for all $1\leq k,l, k', l'\leq n$,
\begin{eqnarray*}
	\mathbb E Y_{kl}&=& \delta_{kl}+S_\nu(E_k,E_l)\\
	\text{ and } \text{Cov}\left(Y_{kl}Y_{k'l'}\right)&=& (-2R +I^2)\left((E_k,E_l),(E_{k'},E_{l'})\right).\label{covid}
\end{eqnarray*}
\item
Assume now that $n=N$. In this case, there is no normal bundle.
Let us define the equivalent of $\mu_{t,\nu}$:
\begin{equation}\label{mutnu2}
\forall t\in \R, x\in M,\ 	\mu_{t}\sim 
N(-tI , -2R +I^2).
\end{equation}
\end{itemize}
The following numbers will quantify the error terms in the main theorems. 
Let
\begin{eqnarray}\label{sigma}
\si(M) &=& \sup_{x\in M}\rho(-2R+I^2),
\end{eqnarray}
where $\rho(P)$ denotes the spectral radius of $P$. 
which are positive by continuity of the terms and compacity of $\overline M$, and where $\|\cdot\|$ denotes the norm associated to the standard metric over $\R^n$.
Define also  
\begin{equation}\label{rho}
\rho(M)= \inf_{x\in  M} |\det (-2R+I^2)|^{1/2}
\end{equation}
and 
\begin{equation}\label{esse}
s(M)= \left\lbrace\begin{array}{ccc}
\underset{(x,\nu) \in S T^\perp M}{\sup}\rho(S_{\nu}) &\text{ if } &n\leq N-1\\
0 &\text{ if } &n=N,
\end{array}\right.
\end{equation}
where $S T^\perp M$ 
is the spherical normal bundle, that is 
\begin{equation}\label{sphericalbundle}
\forall x\in M,\ 
ST_x^\perp M = \{\nu \in T_x^\perp M, \|\nu\|_g = 1\}.
\end{equation}
The following constant will measure the exponential decay 
of the critical points which are non-maxima:
\begin{equation}\label{cem}
\theta^{-1} (M)= \max\left(s^2+\si, (s+1)^2\right).
\end{equation}
Define also
\begin{equation}\label{uzero}
u_0(M)=(1+s(M))\max\left(1,\sigma(M)^{1/2}\right)
\end{equation}
and
\begin{equation}\label{uun}
u_1(M)=\max\left(u_0(M), \frac{1}\theta(N^2+2)\right).
\end{equation}
Notice that $\mu_t$, $u_0$, $u_1$, $s$, $\rho$, $\sigma$ and $\theta$ depend only on the $4$-jet of $e$ on the diagonal of $M\times M$. Indeed, $g$ depends on the $2$-jet, and the curvature and second form depend on the second derivatives of the metric. 

We will need a version of the latter parameters not only for submanifolds, but also for stratified sets. Hence, assume that $M\subset \widetilde M$ is a stratified set of dimension $n$. For any of the parameters $\varphi$ defined above, let
\begin{equation}\label{phi}
\varphi(M)= \sup_{j\in \{0,\cdots, n\}}\varphi(\partial_j M) \text{ or }
\inf_{j\in \{0,\cdots, n\}}\varphi(\partial_j M) 
\end{equation}
depending of the definition of $\varphi$.

\paragraph{A general formula.}
The following lemma provides a general Kac-Rice formula for the average of critical points of given index. For this, let $n\geq 1$ and 
\begin{align}\label{syme} 
\forall 1\leq i\leq n,\ 
\sym(i,n-i,\R)&= \{A\in \sym(n,\R), \Ind(A)=i\},
\end{align}
where $\Ind(A)$ is defined by~(\ref{index}).
 In the sequel, for any topological subspace $M\subset \widetilde M$, $$\partial M=\overline{M}\setminus M.$$
Recall that if $M$ is the interior of a submanifold with boundary, this boundary coincides with the geometrical boundary. If $M$ is the stratum of stratified sets, $\partial M$ can be more intricated.

\begin{lemma}\label{critman}Let $\widetilde M$ a $C^3$ manifold of dimension $N$,  
	$n\in \{0,\cdots, N\}$ be an integer and $M\subset \widetilde M $ be a $C^2$ $n$-dimensional submanifold of $M$, such that  $\partial M$ has finite $(n-1)$-th Hausdorff measure. Let $\widetilde f$ be a centered Gaussian field on $\widetilde M$ satisfying the
	conditions~(\ref{C1}) (regularity),  (\ref{C2}) (non-degneraticity) and~(\ref{C4}) (unit variance), $g$ be the induced metric defined by~(\ref{metric}), $f=\widetilde f_{|M}$, 
	and $i\in \{0,\cdots, n\}$. Then,
	\begin{eqnarray*} \forall u\in \R, \ \mathbb E C_{i}(S_u(M,f))& =&  \frac{1}{\sqrt{2\pi}^{N+1}}
		\int_{(x,\nu)\in T^\perp M} \int^{u}_{-\infty} e^{-\frac{1}2t^2}  e^{-\frac{1}2\|\nu\|_g^2} \\
		&& \int_{Y\in \sym(i,n-i,\R)}|\det Y| d\mu_{t,\nu}(Y)dt d\vol_{g}(x,\nu) ,
	\end{eqnarray*}
	where $d\mu_{t,\nu}$ is defined by~(\ref{mutnu}). If $n=N$, the integral in $\nu$ is removed and $\mu_{t,\nu}$ is replaced by $\mu_t$ defined by (\ref{mutnu2}).
\end{lemma}
\begin{preuve} The Kac-Rice formula given by~\cite[Theorem 12.1.1]{adler} is written for compact manifolds but holds for open manifold whose topological boundary has finite $(n-1)$-Hausdorff measure, as in~\cite[Theorem 11.2.1]{adler}. Now, from the proof of~\cite[Theorem 12.4.2]{adler}, for all $i\in \{0\cdots, n\}$ and every $ u\in \R, $ we get
	\begin{eqnarray*}
		\mathbb E C_{i}(S_u(M,f)) &=& 
		\int_{(x,\nu) \in T M^\perp}
		\mathbb E
		\left(|\det \nabla^2 f_{|T_xM}|  {\bf 1}_{ f(x)\leq u} {\bf 1}_{\Ind (\nabla^2 f_{|T_xM})=i}  \ \Big| \   df_{|T_xM\oplus T_x^\perp M}=(0,\nu^*)\right)\\
		&&p_{d f(x)}(0,\nu^*) d\vol_g (x,\nu),
	\end{eqnarray*}
	where $\det \nabla^2 f_{|T_xM}$ denotes the determinant 
	of the matrix of the bilinear form $\nabla^2 f_{|T_xM}$ in some orthonormal (for $g$) basis
	of $T_xM$, $\nu^* = g(\nu, \cdot)$, and $p_{df(x)}$ denotes the Gaussian density of $df(x)$.  Using the
	independence of $f(x)$ and $df(x)$ induced by 
	the constance of the variance of $f$,  
	the integrand of the integral over $M$ rewrites
	$$ \frac{1}{\sqrt{2\pi}^{N+1}}\int^{u}_{-\infty}\int_{\nu \in T_x^\perp M} e^{-\frac12t^2}  e^{-\frac{1}2\|\nu\|_g^2} \theta_x(\nu,t) d\vol_{g}(\nu)dt,$$
	where
	$$ \theta_x(\nu,t) =
	\mathbb E 
	\left(|\det \nabla^2 f_{|TM}|\,{\bf 1}_{\Ind(\nabla^2f_{|TM}) = i } \Big| \ f(x)=t, df(x)_{|TM\oplus T^\perp M}=(0,\nu^*) \right).
	$$
	By Proposition~\ref{gorille}, 
	$$ \theta_x(\nu,t) =
	\int_{Y\in \sym(i,n-i,\R)}|\det Y| d\mu_{t,\nu}(Y),$$
	where $\mu_{t,\nu}$ is the Gaussian measure 
	with average $-tI + S_\nu$ and variance
	$ (-2R +I^2)$ at $x$ and restricted to $T_xM$. The case $n=N$ is the same, except that the $\nu$ part is absent. 
\end{preuve}
\paragraph{The total number of the critical points.}
Corollary~\ref{critz} below provides an asymptotic equivalent of the average sum of critical points, with a quantitative error bound.
Before writing it, 	under the hypotheses of Lemma~\ref{critman}, we define
\begin{align}\label{ve}
v(M) &= \left\lbrace\begin{array}{ll}
\vol_g  (M) &  \text{ if } n=N\\
\frac{\vol_g(S T^\perp M)}{\vol_{g_0}(\mathbb S^{N-n-1})} 
&\text{ if } n\leq N-1 
\end{array}\right.,
\end{align}
where $g_0$ 
denotes the standard metric on $\R^{N-n}$.

\begin{corollary}\label{critz}
	Under the hypotheses of Lemma~\ref{critman},
	\begin{eqnarray*}
		\forall u\leq -1, \
		\mathbb E C(S_u(M,f))&=& 	\frac{1}{\sqrt{2\pi}^{n+1}}v(M)|u|^{n-1} e^{-\frac12 u^2}+ \epsilon_{n,u},
	\end{eqnarray*}
	where $v(M)$ is defined by~(\ref{ve}) 	and 
	\begin{equation}\label{epsi}
\forall u\leq -1,\ 	
|\epsilon_{n,u}|\leq P_N(\rho, \sigma_+)v(M) |u|^{n-2}e^{-\frac12 u^2}.
	\end{equation}
	Here, $P_N$ is a real polynomial with non-negative coefficients depending only on $n$, and $\rho, \sigma_+$ and $s$ are defined by~(\ref{rho}), (\ref{sigma}) and ~(\ref{esse}).
\end{corollary}
\begin{preuve}Assume first that $n\leq N-1$.  Let $x\in M$, $t\in \R$, $\nu \in T_x^\perp M$ and $\mu_{t,\nu}$ defined by~(\ref{mutnu}).
	Then, 
	$$ \int_{\sym(n,\R)}|\det Y| d\mu_{t,\nu}(Y)=
	\int_{\sym(n,\R)}|\det (tI+Y-S_\nu)| d\widetilde \mu(Y),$$ 
	where $d\widetilde \mu$ is the centered Gaussian measure with covariance
	$ (-2R+I^2)$, see \S~\ref{indu} for the definitions of $R$, $I$ and $S_\nu$. Here, there is an abuse of notation, since $I, S_\nu$ and $R$ are considered through a fixed orthonormal basis or $T_xM$. Now for two matrices $A,B$ of size $ n$ and a non negative real $t$, 
	\begin{equation}\label{detb}|\det (tA+B) - t^ n \det A |\leq  n! \sum_{i=1}^{n} t^{ n-i} \|A\|_\infty^{ n-i}\|B\|_\infty^{i},
	\end{equation}
	where $\|(A_{kl})_{1\leq k,l\leq  n}\|_\infty=\min_{k,l} |A_{kl}|.$
	Moreover, by H\"older's inequality, there exists positive constant $C_ n, C'_ n$ depending only on $ n$ such that for any $i\in \{1, \cdots,  n\}$,
	\begin{eqnarray*}
		\int_{\sym( n,\R)} \|Y-S_\nu\|_\infty^{i}d\widetilde \mu(Y)  &\leq&
		C_ n\left(\int_{\sym( n,\R)} \|Y\|_\infty^{i}d\widetilde \mu (Y) +\|S_\nu\|_\infty^{i}\right)\\
		& \leq & C'_ n( \sigma_+^{i/2}+\|\nu\|_g^i s^i).
	\end{eqnarray*}
	Here, we used that 
	\begin{equation}\label{nunuche}
	\forall \nu\in T^\perp_x M\setminus\{0\},\   S_\nu = \|\nu\|_g S_{\nu/\|\nu\|_g}.
	\end{equation}
	Hence, by~(\ref{detb}) and Lemma~\ref{critman}, there exists 
	a polynomial $P_N$ depending only on $N$ and with non-negative coefficients, such that  for all $u\leq -1$, 
	$$
	\mathbb E C(E_u)- 	\frac{|u|^{ n-1} e^{-\frac12u^2}}{\sqrt{2\pi}^{N+1}}
  \int_{(x,\nu) \in T^\perp M} e^{-\frac12\|\nu\|_g^2} 
	d\vol_{g}(x,\nu)
	$$ is bounded by
	$|u|^{ n-2}e^{-\frac12 u^2} P_N(\rho,\sigma_+)\vol_g	 (ST^\perp M).$
	Now using the coarea formula applied to $v\mapsto \|v\|_g$, 
	\begin{eqnarray*}
		\int_{\nu \in T_x^\perp M} e^{-\frac12\|\nu\|_g^2} 
		d\vol_{g}(\nu)&=&\int_0^{+\infty} t^{N- n-1}e^{-\frac12 t^2}dt \vol_g(S_x T^\perp M)\\
		& = & \sqrt{2\pi}^{N- n} \frac{\vol_g(S_x T^\perp M )}{\vol_{g_0}(\mathbb S^{N- n-1})}.
	\end{eqnarray*}
	For $ n=N$, the normal bundle is the zero space, so that 
	the latter integral must be considered to be equal to 1.
\end{preuve}

\paragraph{Non-maxima critical points.}
We prove now a quantitative version of 
\cite[Theorem 5.2.1]{adler10} for submanifolds, which can be hence implemented immediatly into the general context of stratified sets. It says that in average, the proportion of critical points in $E_u(M,f)$ which are not local maxima decreases exponentially fast with $u$. 
\begin{theorem}\label{critical}
	Under the hypotheses of Lemma~\ref{critman},
	\begin{eqnarray*}
		\forall u\leq -u_0, \ 
		\mathbb E C(S_u(M,f))=\mathbb E C_0(S_u)+ \eta_{ n,u},
	\end{eqnarray*}
	where 
	\begin{equation}\label{eta}
		\forall u\leq -u_0, \ 	|\eta_{ n,u}|\leq  	
	\frac{1}{\rho}Q_N(\sigma^{1/2}, \sigma^{-1/2}, s)v(M)|u|^{N^2}e^{-\frac{1}2u^2(1+\theta)},
	\end{equation}
	where $v(M)$ is defined by~(\ref{ve}), $u_0$ by~(\ref{uzero}), $\sigma$ by~(\ref{sigma}), $\rho$ by~(\ref{rho}), $s$ by~(\ref{esse}), $\theta$ by~(\ref{cem}) and $Q_N$ is a real polynomial depending only on $N$ with non-negative coefficients. 
\end{theorem}
\begin{preuve}
	Lemma~\ref{critman} implies that for all $u\in \R$, 
	\begin{eqnarray*}
		\mathbb E (C(S_u)-C_0(S_u))&= &
		\frac{1}{\sqrt{2\pi}^{N+1}}
		\int_{(x,\nu)\in T^\perp M} \int^{u}_{-\infty} e^{-\frac12t^2}  e^{-\frac12\|\nu\|_g^2} \\
		&& \int_{Y\in \bigcup_{i\geq 1}\sym(i, n-i,\R)}|\det Y| d\mu_{t,\nu}(Y)d\vol_g(x,\nu)dt ,
	\end{eqnarray*}
	where $\mu_{t,\nu}$ denotes the Gaussian measure on $\sym( n,\R)$ given by~(\ref{mutnu}).
	In the sequel, the scalar product on $\sym( n,\R)$
	is $\langle R, S\rangle = \Tr (RS)$.
	Then the integral in $Y$ is bounded by
	\begin{equation}\label{integrale}
	\int_{Z\in \bigcup_{i\geq 1}\sym(i, n-i,\R) +t I-S_\nu}
	|\det (Z+tI-S_\nu)|
	e^{-\frac{1}2\sigma_+^{-1}\|Z\|^2}\frac{dZ}{\sqrt{2\pi}^{ n( n+1)/2} \rho}.
	\end{equation}
	Now if 
	$\rho(S)$ denotes the spectral radius of $S$, then  
	\begin{equation}\label{spec0}
	Z \in \bigcup_{i\geq 1}\sym(i, n-i,\R) +t I-S_\nu\Rightarrow 
	\text{Spec}(Z) \cap \left(-\infty, t +s\|\nu\|_g\right) \neq \emptyset,
	\end{equation}
	where we used~(\ref{nunuche}).
	Moreover by H\"older inequality, there exists a  constant 
	$C_ n>0 $ depending only on $ n$, such that 
	$$ |\det (Z+tI-S_\nu)|\leq  \left(\rho(Z)+|t|+s\|\nu\|_g)^ n\leq C_ n (\rho^ n(Z)+|t|^ n+s^ n\|\nu\|_g^ n\right).
	$$	
	Let  $\lambda_1\leq \cdots \leq \lambda_ n$ be the eigenvalues of $Z$. By (\ref{spec0}),
	\begin{equation}\label{spec}
	Z \in \bigcup_{i\geq 1}\sym(i, n-i,\R) +t I-S_\nu\Rightarrow 
	\lambda_1 \leq t +s\|\nu\|_g.
	\end{equation}
	Writing $Z$ as $$Z={}^tP \text{Diag}(\lambda_1, \cdots, \lambda_ n)P$$ with $P$ orthogonal and using the coarea formula~\cite[Theorem 2.5.2]{anderson}, there exists a  constant $C'_{ n}>0$ such that 
	the integral (\ref{integrale}) is bounded above by
	\begin{eqnarray*}
		\frac{C'_{ n}}{\rho}
		\int_{
			{(\lambda_i)_{2\leq i\leq  n}\in \R^{ n-1}\atop \lambda_1\leq t +s\|\nu\|_g}}
		(\|\lambda\|_{\infty}^ n+|t|^ n+s^ n\|\nu\|_g^ n)
		\|\lambda\|_{\infty}^{ n( n-1)/2}
		e^{-\frac12 \sigma^{-1}\|\lambda\|_2^2}d\lambda,
	\end{eqnarray*}
	where $\lambda=\left( \lambda_1, (\lambda_i)_{2\leq i\leq  n}\right)$
	and $d\lambda$ is the associated Lebesgue measure. 
	After the change of variables $$\mu=\sigma^{-1/2} \lambda,$$ we see that the integral is bounded by
	\begin{eqnarray*}
		\frac{C'_{ n}}{\rho}
		\int_{
			{\mu'\in \R^{ n-1}\atop \mu_1\leq \si^{-1/2}(t +s\|\nu\|_g)}}
		\left(\si^{ n/2}(\|\mu'\|_{\infty}^ n+|\mu_1|^ n)+|t|^ n+s^ n\|\nu\|_g^ n\right)
		\left(\|\mu'\|_{\infty}^{ n( n-1)/2}+|\mu_1|^{ n( n-1)/2}\right)
		\\
		\si^{ n( n-1)/2+ n/2}
		e^{-\frac{1}2\|\mu\|_2^2}d\mu,
	\end{eqnarray*}
	where $\mu=(\mu_1, \mu')$.
	Hence, after integrating in $\mu'$, we see that there exists another constant $C''_ n$, such that the integral is bounded by
	\begin{eqnarray*}
		C''_{ n}			\frac{\si^{ n^2/2}}{		\rho}
		\int_{
			{ \mu_1\leq \si^{-1/2}(t +s\|\nu\|_g)}}
		\left(\si^{ n/2}(1+|\mu_1|^ n)+|t|^ n+s^ n\|\nu\|_g^ n\right)
		(1+|\mu_1|^{ n( n-1)/2})
		e^{-\frac{1}2\mu_1^2}d\mu_1.
	\end{eqnarray*}
	For $s>0$, the latter is bounded by $C'''_{ n}\frac{\si^{ n^2/2}}{\rho}(\sigma^{ n/2}+|t|^ n+s^ nr^ n)(f_1+f_2)(r,t)$, where
	\begin{eqnarray*}
		f_1 &=& {\bf 1}_{\sigma^{-1/2}(t+sr)\geq -1}\\
		f_2 &=& |\sigma^{-1/2}(t+sr)|^{ n( n+1)/2-1}
		e^{-\frac12\sigma^{-1}(t+sr)^2}{\bf 1}_{\sigma^{-1/2}(t+sr)<-1},
	\end{eqnarray*}
	where $r=\|\nu\|_g$ and $C'''_ n$ depends only on $ n$. 
	We used the fact that for any $k\in \Nn$, there exists a constant $c_k$ depending only on $k$, such that 
	$$ \forall u\leq -1, \ \int_{-\infty}^u |x|^k e^{-\frac12 x^2 } dx \leq c_k |u|^{k-1} e^{-\frac12 u^2}.$$
	Consequently, there exists a constant $C_N$ depending
	only on $N$, such that 
	$ \mathbb E (C-C_0)$ is bounded by
	\begin{eqnarray*}
		C_N\vol_g (S T^\perp M)\frac{ \sigma^{ n^2/2}}{\rho} 
		\int^u_{-\infty} \int_0^\infty e^{-\frac12(t^2+r^2)}
		(\si^{ n/2}+|t|^ n+s^ nr^ n)(f_1+f_2)
		r^{N- n-1}dr  dt.
	\end{eqnarray*}
	The double integral splits into two sums $J_1$ and $J_2$, one for $f_1$ and the other for $f_2$. 
	Assume from now on that $$u\leq -u_0=-(1+s)\max(1,\sigma^{1/2}).$$ The bound by the $\sigma^{1/2}$ term implies that 
	$$\forall t\leq u\  \forall r\geq 0, \  \sigma^{-1/2}(t+sr)\geq -1\Rightarrow r\geq -t/(s+1).$$ Changing $t$ into $-t$, there exists $D_N, D'_N>0$ depending only on $N$, such that, using the bound $t\geq 1+s$ and then $u\leq -1$, 
	\begin{eqnarray*}
		J_1& \leq &
		\int_{-u}^\infty 
		\int_{t/(s+1)}^\infty e^{-\frac12(t^2+r^2)}
		(\si^{n/2}+t^n+s^nr^n)r^{N-n-1}
		dr dt\\
		& \leq & 
		D_N \int_{-u}^\infty 
		e^{-\frac{1}2t^2(1+1/(s+1)^2))}
		(\si^{n/2}+t^n+s^n)	(t/(s+1))^{N-2}
		dt \\
		& \leq & 
		D'_N
		(\si^{n/2}+1+s^n)	u^{N+n-3}
		e^{-\frac12u^2(1+\frac{1}{(s+1)^2}) }.
	\end{eqnarray*}
	The second integral $J_2$ satisfies
	\begin{eqnarray*}
		J_2 \leq 
		\int_{-u}^{+\infty} \int_0^{\infty}& e^{-\frac{t^2}2-\frac{r^2}2-\frac{1}2\si^{-1}(t-sr)^2}
		(\si^{n/2}+t^n+s^nr^n)\\
		&	\si^{-(n^2+n)/4}(t+sr)^{n(n+1)/2}
		r^{N-n-1} drdt,
	\end{eqnarray*}
	where we  used $|\sigma^{-1/2}(t-sr)|^{n(n+1)/2-1}\leq |\sigma^{-1/2}(t+sr)|^{n(n+1)/2}$
	in order to simplify the power.
	We write
	$$ t^2+r^2+{\si}^{-1}(t-sr)^2= 
	\left(1+\frac{\si^{-1}}{1+s^2\si^{-1}}\right)t^2+(1+s^2\si^{-1})R^2,$$
	with $R= |r-s\si^{-1} t|.$ 
	Using again H\"older, there exists $E_N$, $E'_N$ and $E''_N$ depending only on  $N$ such that, using $1+s^2\si^{-1}\geq 1$ in the exponential, 
	\begin{eqnarray*}
		J_2& \leq & E_N \si^{-n(n+1)/4}
		\int_{-u}^\infty \int_0^{+\infty} e^{-\frac{t^2}2(1+\frac{\si^{-1}}{1+s^2\si^{-1}})-(1+s^2\si^{-1})\frac{R^2}2}
		(\si^{n/2}+t^n+s^nR^n +s^{2n}\si^{-n} t^n )\\
		&&	\left(
		(sR)^{n(n+1)/2}+((1+s^2\si^{-1}) t)^{n(n+1)/2}\right)
		\left(R^{N-n-1}+(s\si^{-1} t)^{N-n-1}\right) dRdt\\
		&\leq &
		E'_{N} \si^{-n(n+1)/4}
		\int_{-u}^\infty  e^{-\frac{t^2}2(1+\frac{\si^{-1}}{1+s^2\si^{-1}})}
		\left(\si^{n/2}+t^n+s^n +s^{2n}\si^{-n} t^n )\right)\\
		&&	\left(
		s^{n(n+1)/2}+((1+s^2\si^{-1}) t)^{n(n+1)/2}\right)
		\left(1+(s\si^{-1} t)^{N-n-1}\right)dt\\
		&\leq & E''_{N}\si^{-n(n+1)/4}\left(\si^{n/2}+1+s^n +s^{2n}\si^{-n} \right)(s^{n(n+1)/2}+1+(s^2\si^{-1})^{n(n+1)/2})
		\\
		&&
		\left(1+(s\si^{-1})^{N-n-1}\right)|u|^{n(n+1)/2+N-2}
		(1+{\si^{-1}})^{n(n+1)/4+N/2}
		e^{-\frac{u^2}2(1+\frac{\si^{-1}}{1+s^2\si^{-1}})}.
	\end{eqnarray*}
	Summing $J_1$ and $J_2$ implies that for $u\leq -u_0$, 
	\begin{eqnarray*} \mathbb E (C(S_u)-C_0(S_u))&\leq&
		|u|^{N^2}e^{-\frac{u^2}2\left(1+\min(\frac{1}{s^2+\si}, \frac{1}{(s+1)^2})\right)}\vol_g (ST^\perp  M)
		\frac{1}{\rho}P_N(\si^{1/2}, \si^{-1/2}, s) ,
	\end{eqnarray*}
	where $P_N$ is a polynomial depending only on $N$.
\end{preuve}

 \subsection{The main theorem for manifolds}\label{mtm}
We can now prove Theorem~\ref{sph0} for compact manifolds  without boundary, namely Corollary~\ref{coroman}.
In fact, the proof we give is the one of the more precise Theorem~\ref{sph03} 
 below, again for manifolds. 
\begin{preuve}[ of Corollary~\ref{coroman}]
	By Corollary~\ref{coromilnor} and Theorem~\ref{johnny},
	for any $u\in \R$,
	\begin{equation*}
	\sum_{i\geq 1} \mathbb E C_i(S_u(M))\leq \mathbb E C(S_u)-\mathbb E N(S_u)\leq \mathbb E C(S_u)-\mathbb E N_{\mathbb B^n}
	(S_u)\leq
	2	\sum_{i\geq 1} \mathbb E C_i(S_u).
	\end{equation*}

	By Theorem~\ref{critical},
	$$u\leq -u_0, \ \sum_{i\geq 1} \mathbb E C_i(S_u)\leq \eta_{n,u},$$
	where $\eta_{n,u}$ satisfies the bound~(\ref{eta}). 
	Besides, by Corollary~\ref{critz},
	\begin{eqnarray*}
		\forall u\leq -1, \
		\mathbb E C(S_u(M,f))&=& 	\frac{1}{\sqrt{2\pi}^{n+1}}\vol_g(M)|u|^{n-1} e^{-\frac12 u^2}+ \epsilon_{n,u},
	\end{eqnarray*}
	where $\epsilon_{n,u}$ satisfies the bound~(\ref{epsi}). 
	Consequently, 
	\begin{eqnarray*}
		\forall u\leq -u_0, \
		\mathbb E N_{\mathbb B^n} (S_u(M,f))&=& 	\frac{1}{\sqrt{2\pi}^{n+1}}\vol_g(M)|u|^{n-1} e^{-\frac12 u^2}+ \delta_{u},
	\end{eqnarray*}
	where  for $u\leq -u_0$,
	$$|\delta_u|\leq  |\epsilon_{n,u}|+2|\eta_{n,u}|\leq   Q_N(\rho^{-1},\rho, \si^{-1/2}, \si^{1/2},s)\vol_g(M)(|u|^{N^2}e^{-\theta u^2}+ 
	|u|^{n-2})e^{-\frac12u^2},$$
	where $Q_N$ is a real polynomial depending only on $N$ and with non-negative coefficients
	and $\theta $ given by~(\ref{cem}).
	 Hence 
the first assertion of Corollary~\ref{coroman} is proven. 

	We turn now to the second assertion concerning $N_{\mathbb S^{n-1}}(Z_u(M,f)$.
	By Corollary~\ref{pasteque}, 
	For any $u\in \R$, 
	$$| \mathbb E C(S_u)-\mathbb E N(Z_u(M))|\leq 7\sum_{i=1}^n \mathbb E  C_i(S_u),$$
	and the same holds for $N_{\mathbb S^{n-1}}$ instead of $N$.
	We conclude as before.
\end{preuve}

\section{Stratified Morse theory}\label{morsestrat}

We present in this section Morse theory for Whitney stratified set. 

\subsection{Whitney stratified sets} 
A  stratified set is a disjoint union of manifolds satisfying certain gluing conditions. 
\begin{definition}\cite[p. 185]{adler},\cite[pp. 36--37]{goresky}\label{stratified} Let $k\geq 1$ be an integer and $\widetilde M$ be a $C^k$ manifold.
	For  $1\leq \ell \leq k$, a  $C^\ell$ \emph{stratified space} $(M,Z)$ is a pair of a subset $M\subset \widetilde M$  equipped with a 
	locally finite partition $Z$ of $M$ satisfying the following conditions:
	\begin{enumerate}[series=condi]
		\item 
		each piece, or \emph{stratum}, $S \in Z$ is a $C^\ell$ submanifold of $\widetilde M$, without
		boundary and locally closed;
		\item for $R, S \in Z$, if $R \cap \overline{S} \neq \emptyset$, then $R \subset \overline{S}$, and $R$ is said to be \emph{incident} to $S$.
	\end{enumerate}
	A stratified space $(M,Z)$ is said to be \emph{Whitney}
	if it satisfies the further condition:
	\begin{enumerate}[resume=condi]
	\item for any $R, S\in Z$, $R$ incident to $S$, any $x\in R$, any sequences $(x_n)_n \in R^\Nn$ and $(y_n)_n \in S^\Nn$
	with $x_n \to x$ and $y_n \to x$, and if $\varphi : \widetilde M \to \R^N$ is a local chart near $x$, 
	the sequence of line segments $[\varphi(x_n )\varphi(y_n )]\in \R^N$ converges in projective
	space to a line $\ell$  and the sequence of tangent spaces $T_{x_n} S$ 
	converges in the
	Grassmannian to a subspace $\tau \subset T_x \widetilde M$,	 then $d\varphi^{-1} (\ell) \subset \tau$ .
	\end{enumerate}
\end{definition}
One can prove that being Whitney does not depend on the chosen chart of its definition. 
For any stratification $(M,Z)$, the dimension of $M$ is defined by
\begin{equation}\label{dimension}\dim M= \sup_{S\in Z}\dim S.
\end{equation}
The union of strata of dimension $j\in \{0,\cdots, \dim M\}$ 
is denoted by $\partial_j M$, that is 
$$ \partial_j M=\bigcup_{S\in Z, \ \dim S=j} S.$$
\begin{example}\label{trois}
	\begin{itemize}
		\item (Manifolds) If $ M$ is a $C^\ell$ submanifold without (resp. with) boundary of a $C^k$ manifold $\widetilde M$, $k\geq \ell$, then $M$ has a natural structure of $C^\ell$ Whitney stratified space given by $Z=\{M\}$ (resp. $Z= \{M\setminus \partial M, \partial M\}$. 
		\item (Cubes) If $\widetilde M= \R^n$ and $M= [0,1]^n$, then $M$ posses a natural smooth Whitney  stratified structure, with $Z$ being the decomposition into the interior $\mathring{M }$, the interior of its $(n-1)$-faces, etc.  The hypercube $[0,1]^n$ has thus $n+1$ strata, with
		$\partial_n [0,1]^n=]0,1[^n$ and $\partial_{0} [0,1]^n=\{0,1\}^n$.
		\item\label{spiral} (Spirals)	The Whitney condition in Definition~\ref{stratified} is quite subtle: the spiral $$\{e^{-t^2 +it},t>0\}\cup \{0\}\subset \R^2 $$ is a Whitney stratified set, but $\{e^{-t+it}, t\geq 0\}\cup \{0\}$ is not, see~\cite[1.4.8]{pflaum}. 
		Both spirals turn an infinitely number of time around the origin, but the first one is more straight, in a way, than the second one. 
		\item\cite[1.1.11]{pflaum} For a stratified set $L$, the topological cone $$\Cone(L)= ([0,1[\times L)/(\{0\}\times L)$$
		has a natural stratification, where the extremity $[\{0\}\times L]$ is a 0-dimensional
		stratum, and the other strata are $]0,1[\times S$, where $S$ is a stratum of $L$.
	\end{itemize}
\end{example}
See~\cite[p. 187]{adler} for a list of other examples. 


\subsection{Morse functions}

\begin{definition}\label{crit}\cite[p.6, p.52]{goresky}, \cite[p.194]{adler}
	Let $1\leq m\leq \ell\leq k$,  $\widetilde M$ be a $C^k$ manifold, $ \widetilde f : \widetilde M\to \R$ be a $C^k$ map,  $(M, Z)$ be a $C^\ell$ stratified manifold of dimension $n$ and $f:=\widetilde f_{|M}.$
The function $f$ is said to be $C^m$ over 
		$M$ if for $j\in\{0, \cdots, n\}$, 
		$f_{|\partial_j M}$ is $C^m$. 
Moreover, a point $x\in M$ is a \emph{critical point} of $f$  if
		there exists $j$ such that $x$ is a critical point 
		of $f_{|\partial_j M}$. Finally, a critical point $x\in M$ of $f$ is said to be \emph{nondegenerate} if the Hessian of $f_{|\partial_j M}$ in a local chart is non-degenerate at $x$.
\end{definition}
\begin{example}
	Under the hypotheses of Definition~\ref{crit}, any critical point $p\in M$ of $\widetilde f$ in the classical sense is a critical point for $f$ in this wider setting.
If $\dim \widetilde M = \dim M =n$, a point of $\partial_n M$ is critical if and only if  it is critical as a point of $\widetilde M$.  
Any point in $\partial_0 M$ is critical for any function $f$.
\end{example}

\begin{figure}
		\centering
		\includegraphics[width=6cm]{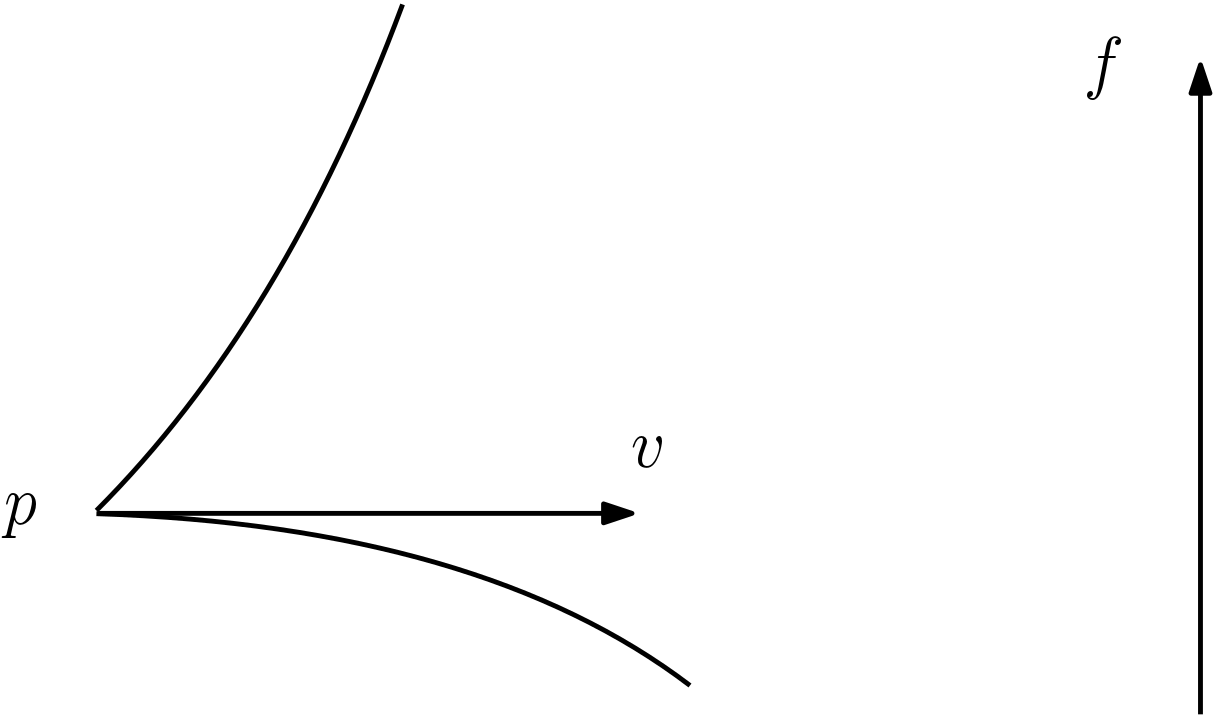}
		\caption{The function $f$ is not Morse over the stratified space $M$, where $\partial_0 M=\{p\}$ and $\partial_1 M $ is the union of the two open drawn branches. Here, $\widetilde M=\R^2$. Indeed $p$ is a critical point of $f$ but $df(p)(v)$ vanishes, where $v$ is a limit of tangents on $M\setminus \{p\}$, hence assertion~\ref{berk}. of Definition~\ref{defmorse} is not satisfied. See also~\cite[pp. 13--14]{goresky}. Rotating a little the figure makes $f$ Morse. 
		}\label{gogo}	
	\end{figure}

\begin{definition}\label{degenerate}~\cite[\S 1.8]{goresky}
	Let $M\subset \widetilde M$ be a $C^1$ Whitney stratified set, 
$x\in M$ and $R$ its stratum. 
\begin{enumerate}
\item  A \emph{generalized tangent space} $T$ at the point $x$ is any subspace of
the form
$$T= \lim_{x_n\to x} T_{x_n} S,$$
where $R\subset \overline S$  and $(x_n)_n$ is a sequence of elements of $S$ converging to $x$ and the limit holds in the Grassmanian of $T\widetilde M$. 
\item A cotangent vector $\nu^*\in T_x^*\widetilde M$  is said to be \emph{degenerate} if 
there exists a generalized tangent space $T\not= T_xR$ such that $\nu^*_{|T}=O.$ Note that in this case, Whitney conditions imply that $\nu^*_{|T_x R}=0$.
\end{enumerate}
\end{definition}
Nondegenerate covectors vanish along $T_x R$ but not in other directions linked to $M$. In Figure~\ref{gogo}, the covector associated to $v$ is degenerate.

\begin{definition}\label{defmorse}\cite[p.52]{goresky}, \cite[Definition 9.3.1]{adler} Let $M\subset \widetilde M$ be a $C^2$ stratified manifold in a $C^2$-manifold, and $\widetilde f : \widetilde M \to \R$ be a $C^2$ map. The restriction $f:= \widetilde f_{|M} : M\to \R$ of $f$ is said to be \emph{a Morse function } if the three following conditions are satisfies:
	\begin{enumerate}
		\item the critical values of $f$ are distinct, that is for any pair of distincts critical points $p,q\in M$, $f(p)\not=f(q)$; 
		\item any critical point  of $f$ in the sense of Definition~\ref{crit} is nondegenerate;
		\item\label{berk} for any critical point $x$ the covector $d\tilde f(x)\in T^*_x\widetilde M$ is nondegenerate in the sense of Definition~\ref{degenerate}.
	\end{enumerate}
\end{definition}
\begin{example}Figure~\ref{gogo} provides a counter example for condition~\ref{berk}. 
	\end{example}
Met $M\subset \widetilde M$ be a $C^1$ stratified set of dimension $n$.
For any subset $E\subset M$, any $i\in \{0,\cdots, n\}$ let $f=\widetilde f_{|M}:M\to \R$ be a Morse function.
As in \S~\ref{morsemm}, 
denote by $\crit_i(E,f)$ the set of critical points 
of index $i$ of $M$ in the sense of Definition~\ref{crit} and belonging to $E$, 
$$C_i(E,f)=\# \crit_i(E,f)\text{ and }
 C(E,f)= \sum_{i=0}^n C_i(E,f).$$
	
\subsection{Change of sojourn sets}\label{dix}
Let $M\subset \widetilde M$ be a $C^1$ Whitney stratified set in a manifold $\widetilde M$ and $x\in M$.  
Let
$D_x\subset \widetilde M$ a small submanifold diffeomorphic to a ball, transverse to the stratum $S$ containing $x$ and such that
$D_x \cap S = \{x\}$. Note that $\dim D_x= \dim \widetilde M-\dim S$.
The associated \emph{normal slice} at $x$ and is defined by~\cite[pp. 7--8]{goresky} 
\begin{equation}\label{slice}
N_x= D_x\cap M.
\end{equation}
\begin{example}
	If $S$ is a neighborhood (in $M$) of $x$, then 
	$N_x=\{x\}$. If $M$ is an $n$-dimensional manifold with boundary and $x\in \partial M$, then $N_x\sim [0,1]$.
	For a canonically stratifed rectangle $M\subset \R^2$, at a corner $x$, $D_x$ is a disc and $N_x$ the quarter of the disc inside $M$. For the example of Figure~\ref{tore-bizarre}, $D_{p_1}$ is a segment and $N_{p_1}=\{p_1\}$, $D_{p_2}$ is a 2-disc and $N_{p_2}$ is the union of three segments glued at $p_2$, $D_{p_3}$ is a 3-ball and $N_{p_3}$ is a union of two cones glued at $p_3$ with a top 2-dimensional disc (the part of the mirror itself). 
\end{example} 
Assume now that $\widetilde M$ and $M$ are $C^2$ and let $g$ be a metric on $\widetilde M$. Let $f=\widetilde f_{|M}$ be a Morse function, where $\widetilde f : \widetilde M\to \R$ is a $C^2$ function.
	For any critical point $p$ in $M$ with critical value $v=f(p)$, we define for $0<\epsilon\ll \delta$: 
	\begin{itemize}
		\item \cite[Definition 3.5.2]{goresky} the \emph{local Morse data}
		\begin{equation}\label{lmd}
		 (A,B)(p,f) = \left(M\cap B(p,\delta)\cap f^{-1}[v-\epsilon , v+\ep], M\cap B(p,\delta)\cap f^{-1}\{v-\epsilon\}\right),
		 \end{equation}
		 where $B(p,\delta)$ is the ball of radius $\delta$ in $\widetilde M$.
\item	\cite[Definition 3.6.1]{goresky} The \emph{normal Morse
		data at }$p$ is the pair of spaces $(A_{N}, B_N)$ 
		$$(A_N,B_N)(p,f)= \left(N_p\cap A, N_p\cap B \right).$$
	We may think of normal
	Morse data at $p$ as Morse data for the restriction of $f$ to the normal slice
	at $p$. 
\item  \cite[Definition 3.6.1]{goresky}
 The \emph{tangential Morse data at} $p$ to be the pair 
	$$(A_T,B_T)(p,f) (A\cap S, B\cap S) .$$  
\end{itemize}
Recall that by Theorem~\ref{milnor}, if $p$ is a critical point of the restriction of $f$ to the stratum $S$ of $p$ and of index $i\in \{0, \cdots, \dim S\}$, then
$$(A_T,B_T)(p,f)  \sim_{hom} (D^i\times D^{\dim S-i}, \partial D^i \times D^{\dim S-i}).$$  
Finally, recall that the product of two topological spaces $A,B\subset X$ is defined by
$$   (A,B)\times (A',B'):=(A\times A', A\times B'  \cup B\times A')\subset X^2.$$
One first important result about stratified Morse theory is the following theorem~\ref{gore3}. 
\begin{theorem}\label{gore3} 
	Let $f$ be a Morse function  on a Whitney stratified space $M$ and $p\in M$ be a critical point of $f$. Then,
	\begin{enumerate} 
		\item \cite[Proposition 3.5.3, Theorem 7.5.1]{goresky}
		the homeomorphic class of the local Morse data $(A,B)(p,f)$ depends does not depend on the choices of Riemannian metric and constants $\epsilon\ll \delta$. 
		\item \cite[Theorem 3.8]{goresky}\label{han} The total space $A_N(p,f)$  of the normal Morse data  is homeomoprhic to the normal slice $N_p$. 
		\item \cite[Theorem 7.5.1]{goresky}	the homeomorphic class of the normal Morse data $(A_N,B_N)(p,f)$ depends only on the differential of $f$ at $p$, and not on the choices of Riemannian metric, the transverse ball $D_x$ and constants $\epsilon\ll \delta$. Moreover, if two differentials are in the same component of the set of nondegenerate covectors, then their associated normal Morse data are also homeomorphic. 
		\end{enumerate}
		\end{theorem}
	Assume that the ambient space $\widetilde M$ is equipped with a metric $g$, and let $M\subset (\widetilde M,g)$ be a stratified set. The \emph{stratified normal bundle}~\cite[p. 195]{adler} of $M$ is defined by
	$$ T^\perp M = \bigcup_{S\in Z} T^\perp S\subset T\widetilde M.$$
	Recall that nondegenerate covectors are defined in Definition~\ref{degenerate}. 
\begin{remark}\label{nunucherie}
Thanks to  Theorem~\ref{gore3}, for any $p\in M$ and any vector $\nu\in T^\perp_p M\subset T_x\widetilde M$, we can associate the normal Morse data  $(A_N,B_N)(p,\nu)$. It is defined by any data of a Morse function having a critical point at $p$ and such that $\nabla \tilde f(p)=\nu$.  The fact that the normal Morse data depends only on the one-jet of $f$ is the reason of the definition (9.2.1) in~\cite{adler}.
\end{remark}

\begin{example}\label{croc}Let  $M\subset \widetilde M$ be a submanifold of dimension $n$ 
with boundary,  equipped with its canonical stratification. Let $x\in \partial M$, $\nu \in T^\perp_x\partial M$ and $n_x\in T_x M$ be a non vanishing outward normal vector to $\partial M$. 
Then 
$$\left\lbrace
\begin{array}{llll}
(A_N,B_N)(p,\nu) &\sim & ([0,1], \emptyset ) &\text{ if } \langle \nu, n_x\rangle <0 \\
(A_N,B_N)(p,\nu) &\sim & ([0,1],\{pt\}) &\text{ if }  \langle \nu, n_x\rangle >0,
\end{array}\right.
$$
where $(A_N,B_N)(p,\nu)$ is defined in Remark~\ref{nunucherie}.
\end{example}

The generalization of Theorem~\ref{milnor} in the stratified setting is the following:
\begin{theorem}\label{Tgoresky}
	Let $\widetilde M$ be a $C^2$ manifold, $M\subset \widetilde M$ be a $C^2$ Whitney stratified space, $\tilde f: M\to \R$ be a $C^2$ map such that $f=\tilde f_{|M}$ is Morse. Then, 
	\begin{enumerate} 
		\item\label{gore1} (Invariance)~\cite[p.6, Theorem SMT Part A]{goresky}  Let $u\leq v\in \R$ be such that $[u,v]$ does not contain any critical value of $f$. Then $$S_u(M,f)\sim_{hom} S_v(M,f)$$ and the intersection of $S_u(M,f)$  with any stratum $S$ is diffeomorphic (up to boundary) to $S_v(M,f)\cap S$. 
		\item\label{boldair}  (Local Morse data) \cite[\S 3.3, Theorem 3.5.4]{goresky} 
		Assume that $p\in M$ is the only critical point in its level set $Z_u(M,f)$, where $f(p)=u$. Let $(A,B)$ be the local Morse data at $p$. Then, for any $\epsilon>0$ small enough, there exists an embedding $g: B\to Z_{u-\epsilon}$ such that 
$$S_{u+\epsilon}(M,f)\sim_{hom}S_{u-\epsilon}(M,f)\cup_g (A,B).$$ 
\item  \label{gore2}	(Factorization of the local Morse data) \cite[p.8, Theorem SMT Part B]{goresky} Under the hypotheses of assertion~\ref{boldair}., 
$$ (A,B)(p,f)=(A_N,B_N)\times (A_T,B_T).$$
	\end{enumerate}
\end{theorem}
\begin{lemma}\label{contractible}
	Under the hypotheses of Theorem~\ref{Tgoresky},
	for a any critical point $x\in M$,
$A$ and $A_N$ are contractible in $M$, where $(A,B)$ and $(A_N,B_N)$ are the local and normal Morse data at $x$.
\end{lemma} 
\begin{preuve}
By~\cite[p. 41]{goresky}, any point of a Whitney stratified set has a neighborhood which is homeomorphic to  the product of a neighborhood of $S$ and a cone, hence contractible in $M$. This implies that $N_x$ is contractible as well. By Theorem~\ref{gore3} assertion~\ref{han}, so is $A_N$. Since $A_T$ is contractible as a product of two balls, so is $A$.
\end{preuve}
\begin{example}\label{exgo}
\begin{itemize}
\item \label{exgo1} If the stratum $S$ of a critical point $p$ is a neighborhood of $p$, $(A_N,B_N)=(\{p\}, \emptyset)$, so that $(A,B)=(A_T,B_T)$ and we recover Theorem~\ref{milnor}. 
\item If $p$ in the stratum $S$ of dimension $j$  is a local minimum for $f_{|S}$  then the local Morse data at $p$ is $(A_N\times D^j, B_N\times D^j)$. 
\item \label{exgo3}
	If $p$ in $S$ of dimension $j$ is a local minimum of $f$ (not only of $f_{|S}$), then $B_N=B=\emptyset$ and the local Morse data is $(A_N\times D^j,\emptyset)$, so that $$S_{u+\epsilon} \sim_{hom} S_{u-\epsilon}\sqcup (A_N\times D^j) .$$ 
	\item In the case where $M$ has a boundary and $p\in \partial M$ is a critical point of index $i$, then 
	by Example~\ref{croc} and Theorem~\ref{Tgoresky},
	$$\left\lbrace
	\begin{array}{llll}
	(A,B) &\sim & ([0,1]\times A_T, [0,1]\times B_T) &\text{ if } \langle \nabla f(p), n_p\rangle <0 \\
	(A,B) &\sim & \left(([0,1]\times A_T), ([0,1]\times B_T) \cup (\{pt\}\times A_T)\right) &\text{ if }  \langle \nabla f(p), n_p\rangle >0,
	\end{array}\right.
	$$
	In the first case, a handle of dimension $i$ is added to $S_{u-\epsilon}$, see also~\cite[Proposition 7.1]{braess}. In the second case, $S_{u-\epsilon}$ is a deformation retract of $S_{u+\epsilon}$~\cite[Proposition 4.1]{braess}, and in fact the are homeomorphic, see Remark~\ref{tark} below. 
	Note that in the latter case $[0,1]\times A_T\sim_{hom}\mathbb B^n$.
	\end{itemize}
\end{example}

\begin{figure}
	\centering
		\includegraphics[width=0.3\textwidth]{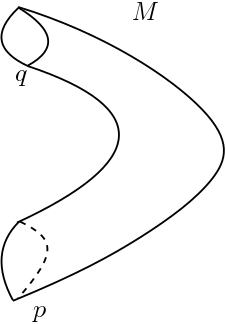}
	\includegraphics[width=0.6\textwidth]{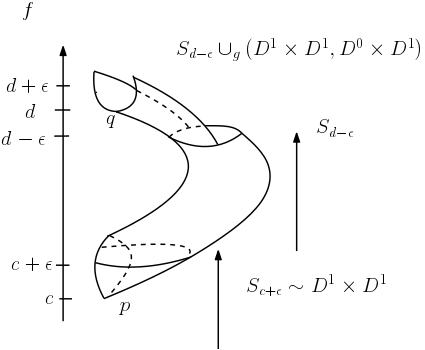}
	\caption{Adding a handle in the case of a manifold with boundary. Here $M$ is a half torus with the canonical stratification of a manifold with boundary. The points $p$ and $q$ are critical for $f$ in the  sense of Definition~\ref{crit}, with vanishing index, so that the tangent Morse data is $(D^0\times D^1, \emptyset)$ in both cases. The normal Morse data at $p$ is $(D^1,\emptyset)$ and $(D^1,D^0)$ at $q$. Theorem~\ref{Tgoresky} asserts that $S_{c+\epsilon}$ is homeomorphic to $S_{u-\epsilon} =\emptyset$ with the handle $(D^1\times D^1, \emptyset)$ attached. In particular,  $S_{u-\epsilon}\sim_{hom} [0,1]\times [0,1]$. Note that it is not so clear how to find an homeomorphism between the local Morse data given by $S_{u+\epsilon}$ and $D^1\times D^1$ seing as the product of the tangent handle and the normal data. This is indeed quite subtle in its full generality and explained in~\cite[8.5]{goresky}. More precisely, this corresponds to the passage from diagram $D_6$ to diagram $D_7$ in~\cite[p.104]{goresky}. At $q$, the handle is $(D^1\times D^1, D^1)$ and $S_{d+\epsilon}\sim_{hom} S_{d-\epsilon}.$ Note that $\langle \nabla f(p),n_p\rangle <0$ and $\langle \nabla f(q),n_q\rangle >0$, see Example~\ref{exgo}. }
	\label{tore-bord}
\end{figure}
\subsection{Morse inequalities}

 We could not find any reference for weak Morse inequalities for stratified sets, but it is quite straightforward from Theorem~\ref{Tgoresky}. If $M$ is a Whitney stratified set, define:
 \begin{align}\label{gamma}
 	\forall i\in \{1, \cdots, n\},\ 
 \gamma_i &= \sup_{\substack{f \text{ Morse}\\ x\in \crit_i(M,f)}} b_{i-1}(B(x,f))\in \R_+\cup \{+\infty\},
 \end{align}
 where $(A,B)(x,f)$ denotes the local Morse data of the critical point $x$. 
\begin{proposition}\label{coq}(Weak Morse inequalities for stratified sets) Under the hypotheses of Theorem~\ref{Tgoresky}, 
	for any non critical value $u\in \R,$
	\begin{eqnarray*} 
			b_0(S_u(M,f))&\leq &C_0(S_u)\\
		\forall i\in \{1, \cdots, n\},\ b_i(S_u(M,f))&\leq& \gamma_i C_i(S_u),
\end{eqnarray*}
where $\gamma_i$ is defined by~(\ref{gamma}).
	\end{proposition}
\begin{example} \label{exgaga}
	\begin{itemize}
\item	In the case where $M$ is a $n$-manifold without boundary and $p$ is a critical point of $f$ with index $i\geq 1$, $B(p,f)\sim \partial D^i\times D^{n-i}$ so that 
$$\forall i\geq 2,\  \gamma_{i} =1$$ and we recover the classical weak Morse inequality given by Theorem~\ref{milnor} for $i\not=1$. For $i=1$, our estimate has a superfluous factor 2. 
\item In the case where $M$ has a boundary and $p\in \partial M$ is a critical point of index $i$, then by Example~\ref{exgo}, $b_{i-1} (B(p,f))=b_{i-1}(\partial D^i)\leq 2$  if $\langle \nabla f(p) ,n_x\rangle <0$; if the latter  is positive, then the Betti numbers of the sojour set does not change. 
\end{itemize}
	\end{example}
Proposition~\ref{coq} is a consequence of the following:
\begin{lemma}\label{relative} Under the hypotheses of Theorem~\ref{Tgoresky}, Let  $i\in \{0, \cdots, n\}$,
and $p$ be a critical point of $f$ with critical value $u\in \R$ and
local Morse data $(A, B)$. Then, for $\epsilon>0$ small enough,
\begin{eqnarray*}
	\dim H_0(S_{u+\epsilon}, S_{u-\epsilon})&\leq &1\\
	H_1(S_{u+\epsilon}, S_{u-\epsilon}) &\hookrightarrow &H_0(B)\\ \text{ and }
 \forall i\geq 2, H_i(S_{u+\epsilon}, S_{u-\epsilon}) 
&\simeq& H_{i-1}(B),  
\end{eqnarray*}
	\end{lemma}
\begin{preuve}
	By the snake lemma, for $i\geq 1$,
	$$\cdots \to H_i(A)\to H_i(A,B)\to H_{i-1}(B)\to 
	H_{i-1}(A)\to \cdots$$
	and
	for $i=0$,
$$\cdots \to H_0(B)\to H_0(A)\to H_0(A,B)\to 0.$$
By Lemma~\ref{contractible}, $A$ is contractible, so that 
$$\cdots \to H_i(\{p\})\to H_i(A,B)\to H_{i-1}(B)\to 
H_{i-1}(\{p\})\to \cdots$$	
hence for $i\geq 2$, 
	$H_i(A,B)\simeq H_{i-1}(B)$,
	for $i=1$, $H_1(A, B)\hookrightarrow H_0(B)$
	and for $i=0$, $H_0(\{p\}) \twoheadrightarrow H_0(A, B)$. 
	 Finally, the excision theorem implies that 
	 $$ H_*(S_{u+\epsilon}, S_{u-\epsilon})\sim H_*(A,B),$$ hence the result.
	\end{preuve}
\begin{preuve}[ of Corollary~\ref{coq}] Let $ i\in\{0,\cdots, n\}.$  
By~\cite[Lemma 5.1]{milnor} applied to $S(A,B)=\dim H_i(A,B)$, for any $\epsilon >0$ small enough, $$\dim H_i(S_u)=\dim H_i(S_u,S_{-\infty})
\leq \sum_{p\in \crit_i(S_u,f)} \dim H_i(S_{f(p)+\epsilon},S_{ f(p)-\epsilon}).$$
By Lemma~\ref{relative}, we obtain the result. 
	\end{preuve}
We provide now the equivalent of Corollary~\ref{coromilnor}. For this, define
\begin{align}\label{beta0}
\beta_0 (M)&= \sup_{\substack{f \text{ Morse}\\ x\in \crit(M,f)}} b_{0}(B(x,f))\in \R_+\cup \{+\infty\}.
\end{align}
\begin{corollary}\label{hibou}
	Under the hypotheses of Theorem~\ref{Tgoresky}, for any non-critical $u\in \R$, 
\begin{equation*}
C_0(S_u(\partial_n M))-  \beta_0(M)\left(C(S_u(M))-C_0(S_u(\partial_n M))\right) \leq N_{\mathbb B^n} (S_u(M))\leq N(S_u(M))\leq C_0(M),
\end{equation*}
	where $\beta_0$ is defined by~(\ref{beta0}).
\end{corollary}
\begin{preuve} Let $x$ be a critical point of $f$ in the sense of Definition~\ref{critical} with $f(x)=u$. Then, by Theorem~\ref{Tgoresky},
	for $0<\epsilon\ll \delta $ small enough,
	$$ S_{u+\epsilon}\sim_{hom} S_{u-\epsilon}\cup_g (A,B)(x,f). $$
	In particular, 
	at most $b_0(B(x,f))$ components of $S_{u-\epsilon}$ can be 
changed. This implies 
	$$N_{\mathbb B^n}(S_{u+\epsilon})\geq N_{\mathbb B^n}(S_{u-\epsilon}) - b_0( B(x,f)), $$ 
	and the same holds with $N$ instead of $N$.
	If $p$ is a critical point with vanishing index in $\partial_n B$, 
	then by Example~\ref{exgo3} and Theorem~\ref{milnor} assertion~\ref{nball}, 
	$$N_{\mathbb B^n}(S_{u+\epsilon})= N_{\mathbb B^n}(S_{u-\epsilon}) +1. $$
	This equality and the former inequality imply the first inequality of the corollary. The second one is trivial, and the last one is due to the first assertion of Corollary~\ref{coq}. 
	\end{preuve}
The next lemma compares two global measures of Morse complexity by a geometrical local one:
\begin{lemma}\label{cracboum}Let $M\subset (\widetilde M,g)$ be a Whitney stratified set in a Riemannian manifold. Then $\beta_0$ and $(\gamma_i)_i$ satisfy the following bounds:
\begin{align} \beta_0 &\leq 2+ \sup_{(x,\nu)\in T^\perp M}
		b_{0} (B_N(x,\nu))\\
	\forall i\in \{1, \cdots, n\}, \  \gamma_i &\leq 2+ \sup_{(x,\nu)\in T^\perp M}
(2	b_{i-1}+b_{i-2}) (B_N(x,\nu)),
	\end{align}
	where $B_N(x,\nu)$ is defined in Remark~\ref{nunucherie}.
\end{lemma}
\begin{preuve}
	Let $x$ be a critical point of $f$. By Theorem~\ref{Tgoresky}, 
	$$B(x,f)\sim_{hom}( A_N\times B_T) \cup (B_N\times A_T).$$	
	Furthermore, if the index of $x$ equal $i\in \{1, \cdots, n\}$,
	$$(	A_N\times B_T) \cap (B_N\times A_T)=
	B_N\times (\partial D^i \times D^{n-i}),
	$$
	which retracts onto $B_N \times \mathbb S^{i-1}$. 
	Now, by Mayer-Vitoris,
	$$
	\cdots \to H_{i-1}(\mathbb S^{i-1})\oplus H_{i-1}(B_N)
	\to H_{i-1} (B) \to H_{i-2}(B_N\times \mathbb S^{i-1})\to \cdots$$
	Hence, $b_{i-1}(B)\leq \dim H_{i-1}(\mathbb S^{i-1})+\dim H_{i-1}(B_N)+\dim H_{i-2}(\mathbb S^{i-1}\times B_N).$
	By the K\"unneth formula, $$b_{i-2}(\mathbb S^{i-1}\times B_N)\leq b_{i-2}(B_N)$$
	so that $b_{i-1}(B)\leq 2+(2b_{i-1}+b_{i-2})(B_N).$
	
	We also have 
	$$\cdots \to 
	H_{0}(\mathbb S^{i-1})\oplus H_{0}(B_N)
	\to H_{0} (B) \to 0,$$
	which implies $ b_0(B) \leq 2+b_0(B_N).$
\end{preuve}

	\subsection{Change of nodal sets}
		In this paragraph we prove a generalization of Proposition~\ref{morseZ} in the case of manifolds with boundary.
	\begin{figure}
		\centering
		\includegraphics[width=0.8\textwidth]{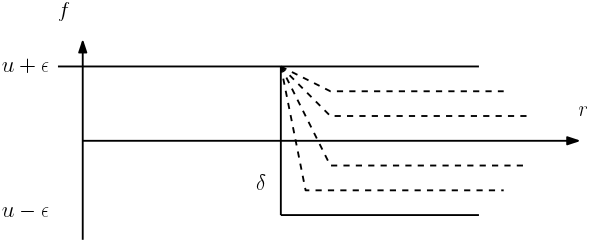}
		\caption{The figure~\cite[p. 96]{goresky}. Here, $r$ is the distance to the critical point $x$ and $f(x)=u$. It shows that $S_{u+\epsilon}$ (the subspace above the horiztonal upper line) is isotopic to $S_{u-\epsilon}\cup_g (A,B)(x,f).$ The dashed paths show the isotopy between the definition of the two subsets. The same argument proves that $Z_{u+\epsilon}$ is homeomorphic to the union of $Z_{u-\epsilon}\setminus B(0,\delta)$ (the upper segment till $\delta$) with $A\cap S(0,\delta)$ (the vertical segment) and with $Z_{u+\epsilon}\cap B(0,\delta)$ (the lower horizontal semi-line).}\label{isot}
	\end{figure}
	We will need a lemma of \cite{goresky} which can be viewed as a nodal version of  Theorem~\ref{Tgoresky} assertion~\ref{boldair}:
	\begin{lemma}\label{toto} 	Let $M\subset \widetilde M$ be a $C^2$ Whitney stratified subset, $f : M\to \R$ be a Morse function and $x\in M$ be a critical point of $f$ with $f(x)=u$. Then, for any $0<\epsilon\ll\delta$ small enough, 
		$$ Z_{u+\epsilon}\sim_{homeo} \left(Z_{u-\epsilon}\setminus B(x,\delta)\right)\bigcup \left((f^{-1}(u+\epsilon)\cap B(x,\delta))\cup
		(f^{-1}[u-\epsilon, u+\epsilon]\cap S(x,\delta))\right),$$
		where $S(0,\delta)$ is the sphere of radius $\delta$ centered on $x$.
		\end{lemma}
	Figure~\ref{isot} shows graphically the proof. 
	\begin{preuve}[ of Lemma~\ref{toto}]
		The proof is a consequence of \cite[\S 7.6]{goresky}. There, it is proven assertion~\ref{boldair}. of Theorem~\ref{Tgoresky},  that is $$S_{u+\epsilon} \sim_{hom} S_{u-\epsilon} \cup_g (A,B),$$ 
		where $(A,B)$ is the local Morse data given by~(\ref{lmd}). But more is proven. From the proof we see that $Z_{u+\epsilon}$ is homeomorphic (even isotopic) to the union of $Z_{u+\epsilon}\cap B(0,\delta)$ (the upper horizontal segment in Figure~\ref{isot}) with $A\cap S(0,\delta)$ (the vertical segment) and with $Z_{u-\epsilon}\setminus B(0,\delta)$ (the lower horizontal semi-line). Hence, the result. 
		\end{preuve}
In Lemma~\ref{cracboum}, using Remark~\ref{nunucherie}, we bounded parameters involving all Morse functions by ones involving only $B(x,\nu)$ for vectors orthogonal vectors $\nu$, that is replacing the general Morse functions by a family of finite dimension $(f_\nu)_{\nu\in T_x^\perp M}$. For nodal sets, would like to do the same, that is replacing the family of Morse functions for a given $x\in M$ by a family of finite dimension. For this, 
let $x\in M$ belonging to the stratum $S$ of dimension $j$.
	Choose $(x_i)_{i\in \{1,\cdots, j\}}$ be local coordinates of $S$ near $x$ and $(y_i)_{i\in \{j+1, \cdots n\}}$ be orthogonal coordinates in $\widetilde M$.
	For $i\in \{0,\cdots, j\}$ and any $\nu\in T^\perp_x M$, 
	let 
	\begin{equation}\label{finu}
	f_{i,\nu}=  -x_1^2-\cdots -x_i^2 +x_{i+1}^2 +\cdots
	+x_j^2 +\nu^*,
	\end{equation}
	where $\forall y\in \{j+1, \cdots, N\}, \ \nu^*(y)=g(y,\nu).$
	\begin{lemma}\label{zedzed}	Let $M\subset \widetilde M$ be a $C^2$ Whitney stratified subset, $f : M\to \R$ be a Morse function, $x\in M$ be a critical point of $f$ of index $i$, and $\nu=\nabla f(x).$ Then, 
		$$ Z_{u+\epsilon}(M,f) \sim_{hom} 
\left(		Z_{u-\epsilon}(M,f)\setminus B(0,\delta) \right)\bigcup \left((f_{i,\nu}^{-1}(\epsilon)\cap B(0,\delta))\cup
		(f_{i,\nu}^{-1}[-\epsilon, \epsilon]\cap S(0,\delta))\right),$$
		where $f_{i,\nu}$ is defined by~(\ref{finu}).
		\end{lemma} 
	\begin{preuve}
		This is a consequence of \cite[Theorem 7.4.1]{goresky}, which asserts that if two Morse functions are isotopic among Morse functions with a unique non-degenerate critical point (here $x$), then their local Morse data are homeomorphic. Here the two functions are $f$ and $f_{i,\nu}$, and it is immediate to check that they satisfy the latter condition. 
	\end{preuve}
Now, let 
	\begin{align}\label{zeta}
	\zeta_0(M):= \sup_{(x,\nu)\in T^\perp M, i\in \Nn} & b_0\left(
	(f_{i,\nu}^{-1}(\epsilon)\cap B(0,\delta))\cup (f_{i,\nu}^{-1}[-\epsilon,\epsilon]\cap S(0,\delta)\right)
	\\
	&  +b_0(f_{i,\nu}^{-1}(-\epsilon)\cap S(0,\delta))+b_0(B(x,\nu))\in \R_+\cup \{+\infty\}.\nonumber
	\end{align}
	Here $0<\epsilon\ll \delta$ are small enough constants depending on $x$. 
By \S~\ref{dix}, $\zeta_0$ does not depend on them. The reason of this parameter is given by the following proposition which generalizes Proposition~\ref{morseZ}.
\begin{proposition}\label{sakharov} Let $M\subset \widetilde M$ be a Whitney stratified set. Let $f : M\to \R$ be a Morse function in the sense of Definition~\ref{defmorse}.
	\begin{enumerate}
		\item	(Invariance)	For any pair of reals $u<v$ such that $f$ has no critical value in $[u,v]$,  $Z_u(M,f) $ is homeomorphic to $Z_v$. 
		\item\label{taka2} (Change at a critical point) For any $u\in \R$, if $p$ is a unique  critical point, in the sense
		of Definition~\ref{crit},  in $Z_u(M,f)$, such that $u=f(p)$, then 
		for $\epsilon$ positive and small enough,
	$$ |b_0(Z_{u+\epsilon}) -b_0(Z_{u-\epsilon})|\leq  \zeta_0,$$
	where $\zeta_0$ is defined by~(\ref{zeta}).
	\end{enumerate}
		\end{proposition} 
\begin{preuve}
	The first point is a consequence of Theorem~\ref{Tgoresky} assertion~\ref{homot}. 
The second assertion follows the global lines of the one of Proposition~\ref{morseZ}. In the sequel, all the subsets involved should be enlarged a little in order  to fit the conditions of Mayer-Vitoris. Since we already wrote a similar proof for Proposition~\ref{morseZ}, we prefer to keep the subsets.  Let
$$ Z_-=Z_{u-\epsilon}\setminus B(0,\delta).$$
Then, by Lemma~\ref{zedzed}, $$Z_{u+\epsilon} = Z_-
 \cup U_{i,\nu},$$
 where $$U_{i,\nu}=\left((f_{i,\nu}^{-1}(\epsilon)\cap B(0,\delta))\cup
(f_{i,\nu}^{-1}[-\epsilon, \epsilon]\cap S(0,\delta))\right),$$
and $\nu=\nabla f(x).$
Then, $$U_{i,\nu}\cap Z_-=V_{i,\nu},$$
where $$ V_{i,\nu}=f_{i,\nu}^{-1}(-\epsilon)\cap S(0,\delta),$$
so that by Mayer-Vitoris, 
$$ H_0(V_{i,\nu})\overset{\alpha}{\to}  H_0(Z_-)\oplus H_0(U_{i,\nu})\to H_0(Z_{u+\epsilon})\to 0,$$
so that $$b_0(Z_{u+\epsilon}) = b_0(U_{i,\nu})+b_0(Z_-)- \rank(\alpha)
$$
In order to bound estimate $b_0(Z_-)$, we note that 
$$Z_{u-\epsilon}= Z_-\cup \left(f_{i,\nu}^{-1}(-\epsilon)\cap B(0,\delta)\right).$$
Recall that $f_{i,\nu} (-\epsilon)\cap B(0,\delta)=B(x,\nu)$
and note that $B(x,\nu)\cap Z_-=V_{i,\nu}.$
Hence, again by Mayer-Vitoris,
$$  H_0(V_{i,\nu})\overset{\beta}{\to}  H_0(B)\oplus H_0(Z_-)\to H_0(Z_{u-\epsilon})\to 0,$$
so that $$b_0(Z_{-}) = b_0(Z_{u-\epsilon})-b_0(B)+ \rank(\beta)
$$
Finally
$$|b_0(Z_{u+\epsilon})-b_0(Z_{u-\epsilon})|\leq b_0(U_{i,\nu}) +2b_0(V_{i,\nu})+b_0(B(x,\nu))\leq \zeta_0,$$
hence the result.
	\end{preuve} 
The next corollary is the equivalent of Corollary~\ref{pasteque}.
\begin{corollary}\label{dauphin}Under the hypotheses of Theorem~\ref{Tgoresky}, 
	$$ \forall u\in \R, \ |N(Z_u(M,f))- C_0(\partial_n M)|\leq \zeta_0(C(S_u)-C_0(S_u(\partial_n M)),$$
	where $\zeta_0$ is defined by~(\ref{zeta}). The same holds for $N_{\mathbb S^{n-1}}$ instead of $N$.
\end{corollary}
\begin{preuve}By Example~\ref{exgo} and By Proposition~\ref{morseZ} assertion~\ref{taka0}, any element of $C_0(\partial_n M)$ creates a new component to $Z_u(M,f)$ diffeomorphic to $\mathbb B^{n-1}
	$. By Proposition~\ref{sakharov}, any other type of critical point can modify the topology of at moste $B(x,f)$ components of $Z_{u-\epsilon}$, and can create at most $\zeta_0$ new components. 
\end{preuve}

We finish this paragraphe proving that for manifolds with boudary, $\zeta_0$ is finite. 
	\begin{proposition}\label{bobord}
	Let $n\in \Nn^*$. There exists $C_n\geq 0$ such that for any 
	$C^2$ manifold $M$  of dimension $n$ and with boundary,  
	$\zeta_0(M)\leq C_n$. 
\end{proposition}
We will need for this the following lemma:
\begin{lemma}\label{coord}\cite[Lemma 3.2]{jankowski1972functions} Let $M$ be a $C^2$ $n$-manifold with boundary, $f: M\to \R$ be a Morse function in the sense of Definition~\ref{defmorse}, and $x\in \partial M$ be a critical point of $f$, of index $i\in \{0, \cdots, n-1\}$.  Then, there exists a neighborhood $U$ of $x$ in $M$ and a coordinate system $(x_i)_{i\in \{1,\cdots, n\} } $ over $U$ such that 
	\begin{itemize}
		\item $\forall p\in U$, $f(p)=x_n(p)$;
		\item $\partial M\cap U =\{		x\in \R^n, \ x_n=-x_1^2 -\cdots -x^2_i +x_{i+1}^2 +\cdots + x_{n-1}^2\}.$
	\end{itemize}			
\end{lemma}
\begin{preuve}[ of Proposition~\ref{bobord}] 
	By Lemma~\ref{coord}, there exists a coordinate system near $x$ such that $f=x_n$ and $\partial M\cap U=\{(x_i)_{1\leq i\leq n}\in x(U), \ x_n=q_i(x_1, \cdots, x_{n-1})\}$, where $q$ is a quadratic polynomial depending only on the index. Since inside this ball $f$ is algebraic, 
	all the subsets defining $\zeta_0$ in~(\ref{zeta}) are semialgebraic (here we use the standard metric on $\R^n$, so that the ball and spheres are algebraic), hence by~\cite{law1965ensembles}, their number of components are finite. 
\end{preuve}

\section{Gaussian fields over Whitney stratified sets}\label{proof}

\paragraph{The conditions for the stratified set.}
In order to apply our results about critical points of random functions to stratified sets and Morse theory to random functions, we will need mild conditions for the stratified set. 
		Let $M\subset (\widetilde M,g)$  be Whitney stratified set in a Riemannian manifold. We sum up all the conditions we need in the article.
	\begin{enumerate}[resume=condition]
\item \label{hausstrat} (gentle boundaries) For any stratum $S$ of dimension $j$, the $(j-1)$-Hausdorff measure of $\partial S$ is finite. 
	\item\label{grass}(very gentle boundaries)\cite[Definition 9.2.1]{adler}   If $S$ is a stratum of $M$, and $p\in \partial S$, then the set $\mathcal T_p=\{\lim_{p_n\to p} T_{p_n} S\}$ of generalized tangent spaces coming from $S$ (see Definition~\ref{degenerate}) has finite Hausdorff dimension less or equal to $\dim S-1$ in the appropriate Grassmannian.
	\item\label{bettistrat0} (mild local connectivity)
$\displaystyle 	\sup_{(x,\nu)\in T^\perp M} b_0(B_N(x,\nu))$ is finite.
		\item\label{bettistrat} (mild local homology)
$\displaystyle \sup_{i\in \{0, \cdots, n\}}		\sup_{(x,\nu)\in T^\perp M} b_{i}(B_N(x,\nu))$ is finite.
	\item~\label{zetastrat} (mild nodal topology) $\zeta_0$ defined by~(\ref{zeta}) is finite. 
	\end{enumerate}
\begin{remark}
Condition~(\ref{hausstrat}) is needed for the Kac-Rice formula, see~\cite[Theorem 11.2.1]{adler}, hence is ubiquitous as far as random critical points are involved. Condition~(\ref{grass}) implies condition~(\ref{hausstrat}) above, and is needed only in the refinement Theorem~\ref{sph2} and Theorem~\ref{euler} from~\cite{adler}. 
Condition~(\ref{bettistrat0}) is needed for Theorem~\ref{sph0} and its quantitative version~\ref{sph03} below.
Condition~(\ref{zetastrat}) is only needed for the nodal versions of the theorems. 
Condition~(\ref{bettistrat}) is only needed for Theorem~\ref{corobetti}. 
\end{remark}

\begin{example}
For a manifold $M$ with boundary and $p\in \partial M$, $\mathcal T_p= \{T_pM\}.$
The spiral given by Example~\ref{trois} does not satisfy condition~(\ref{grass}) but satisfies condition~(\ref{hausstrat}). 
\end{example}

\subsection{The main theorem for stratified sets}\label{mwb}
The following result is the quantitative version of Theorem~\ref{sph0}:
\begin{theorem}\label{sph03}Let $1\leq n\leq N$ be integers,
 $\widetilde M$ be a $C^3$ manifold of dimension $N$,  $M\subset \widetilde M$ be a compact dimension $n$ $C^2$ Whitney stratified set satisfying conditions~(\ref{hausstrat}) (gentle boundaries) and~(\ref{bettistrat0}) (mild local connectivity), $ \widetilde f: \widetilde M \to \R$ be a random centered Gaussian field satisfying conditions (\ref{C1}) (regularity), (\ref{C2}) (non-degeneraticity) and (\ref{C4}) (constant variance), and $f=\widetilde f_{|M}$. Then there exists a polynomial $Q_N$ depending only on $N$ and with non-negative coefficients, such that
	$$ 	\forall u\geq u_1, \ 
	\mathbb E N_{\mathbb B^n}(E_u(M,f)) = \frac{1}{\sqrt{2\pi}^{n+1}} v(\partial_n M)
	u^{n-1}e^{-\frac{1}2u^2} +\delta_{u},$$
	where $v(M)$ is defined by~(\ref{ve})
and	
\begin{equation}\label{deltab}
 	\forall u\geq u_1, \ |\delta_u|\leq  (1+\beta_0) Q_N(\rho^{-1},\rho,\si^{-1/2}, \si^{1/2},s)\max_{1\leq j\leq n}v(\partial_jM)
	u^{n-2}e^{-\frac12u^2}.
	\end{equation}
	The result holds for $N(E_u)$ instead of $N_{\mathbb B^n}$. 
	Here, $\beta_0$, $u_1$, $\sigma$, $\rho$, $s$  and $\theta$ are positive constants depending only on $M$ and $e$ given respectively by~(\ref{beta0}), ~(\ref{uun}), ~(\ref{sigma}), ~(\ref{rho}), (\ref{esse}) and (\ref{cem}). Besides, the volume is computed with respect to the restriction to the $j$-stratum $\partial_j M$ of the metric $g$~(\ref{metric}).
	
	If $M$ satisfies the further condition~(\ref{zetastrat}) (mild local nodal topology), then the same estimate for  $N_{\mathbb S^{n-1}}(Z_u)$ and $N(Z_u)$ instead of $N_{\mathbb B^n}(E_u)$, after
	changing the polynomial $Q_N$ and $\beta_0$ into $\zeta_0$ given by~(\ref{zeta}).
\end{theorem}
 We will need a theorem which asserts that 
under simple hypotheses, a Gaussian random  field is almost surely Morse. 
\begin{theorem}~\cite[Corollary 11.3.5]{adler}\label{johnny}
	Let $\widetilde M$ be a $C^3$ manifold and $M$ be a compact $C^2$ Whitney stratified set, $\widetilde f : \widetilde M \to \R$ a Gaussian random  centered field  and $f=\widetilde f_{|M}$. Then, if $f$ satisfies conditions~(\ref{C1}) (regularity) and (\ref{C2}) (non-degeneracity), then $f$ is almost surely Morse over $M$ in the sense of Definition~\ref{defmorse}.
\end{theorem}
\begin{preuve}Corollary 11.3.5 in~\cite{adler} implies the result if every stratum has a countable atlas. Since $M$ is a Whitney stratified set, then locally there is only a finite number of strata. Moreover, for any stratum $S$, there is an exhaustion of $S$ by a sequence of compacts of $S$, which are all covered by a finite number of charts, hence the result.
	\end{preuve}

\begin{preuve}[ of Theorem~\ref{sph03}]\label{proofstrat}
	The proof is almost the same as in the case of a closed manifold. 
	By Corollary~\ref{hibou} and Theorem~\ref{johnny}, 
	\begin{equation*}
	\mathbb E C_0(S_u(\partial_n M))- \beta_0\mathbb E  \left(C(S_u(M))-C_0(S_u(\partial_n M))\right) \leq  \mathbb E N_{\mathbb B^n} (S_u( M))\leq \mathbb E N(S_u)\leq \mathbb E C_0(S_u).
	\end{equation*}
This implies that $$\max\left(|\mathbb E N(S_u(M))- \mathbb E C(S_u(\partial_n M))|, |\mathbb E N_{\mathbb B^n}(S_u( M))- \mathbb E C(S_u(\partial_n M))|\right)$$
is bounded by
\begin{eqnarray*}
|\mathbb E C_0(S_u(M))-\mathbb E  C(S_u(\partial_n M))|
&+&|\mathbb E C_0(S_u(\partial_n M))-\mathbb E  C(S_u(\partial_n M))|\\
&+&\beta_0|\mathbb E  C(S_u(M))-\mathbb E  C_0(S_u(\partial_nM))|,
\end{eqnarray*}
where $\beta_0$ is defined by~(\ref{beta0}). By Condition~(\ref{bettistrat0}) and Lemma~\ref{cracboum}, $\beta_0$ is finite. 
The first term is bounded by
$$ |\mathbb E C_0(S_u(M))-\mathbb E  C(S_u(M))|+|\mathbb E C(S_u(M))-\mathbb E  C(S_u(\partial_n M))|,$$
where the third term is bounded by
$$
\beta_0(|\mathbb E C(S_u(M))-\mathbb E  C(S_u(\partial_n M))|+|\mathbb E C(S_u(\partial_n M))-\mathbb E  C_0(S_u(\partial_n M))|).
$$
By Corollary~\ref{critical}, for $u\leq - u_0(M)$, the terms involving a difference $C-C_0$ for $M$ or $\partial_n M$
are bounded by $\sum_{j} |\eta_{j,u}|$, 
where $\eta_{j,u}$ satisfies the bound~(\ref{eta}).
This sum is bounded by
$$ 
\frac{1}{\rho}R_N(\sigma^{1/2}, \sigma^{-1/2}, s)\sum_{j\leq n}v(\partial_j M  )|u|^{N^2}e^{-\frac{1}2u^2(1+\theta)},$$
where $R_N$ is a real polynomial with non-negative coefficients and depending only on $N$, and $\sigma$, $\rho$, $s$ and $\theta$ defined by (\ref{sigma}), (\ref{rho}),  (\ref{esse}), (\ref{cem}), ~(\ref{psi}) and $v(M)$ by (\ref{ve}).
By Corollary~\ref{critz}, the terms involving
a difference $C(S_u(M))-C(S_u(\partial_n M))$ are bounded by
 for $u\leq -1$ by 
 $$   \sum_{j\leq n-1}\left( \frac{v(\partial_j M)}{\sqrt{2\pi}^{j+1}} |u|^{j-1}e^{-u^2/2} +\epsilon_{j,u}\right),$$
 where $\epsilon_{j,u}$ satisfies the bound~(\ref{epsi}).
 This sum is bounded by
$$
Q_N(\rho,\sigma)\max_{j\leq n-1}v(\partial_j M)|u|^{n-2}e^{-\frac12 u^2},$$ where $Q_N$ is a polynomial depending only on $N$.  
 Now using that $\forall x\geq 1, \ln x \geq x$,
 	\begin{equation}
 	\forall u\in \R, \ |u|\geq \max \left(1, \frac{1}\theta(N^2-n+2)\right) \Rightarrow |u|^{N^2}e^{-\theta u^2}\leq  
 	|u|^{n-2},
 	\end{equation} the two former bounds imply the first part of Theorem~\ref{sph03}. 
	
	We turn now to the second assertion concerning $N_{\mathbb S^{n-1}}(Z_u(M,f))$.
By Corollary~\ref{dauphin}, 
	$$  \ |\mathbb E b_0(Z_u(M,f))- \mathbb E C_0(\partial_n M)|\leq \zeta_0|\mathbb E (C(M)-C_0(\partial_n M)| ,$$
where $\zeta_0$ is defined by~(\ref{zeta}). By Condition~(\ref{zetastrat}) and Lemma~\ref{cracboum}, $\zeta_0$ is finite.
The rest of the proof is the same as above. 
\end{preuve}

	\begin{preuve}[ of Corollary~\ref{coroman}] The proof for a closed manifold has been done in \S~\ref{mtm}. If $M$ is a $C^2$ compact manifold with boundary, it satisfies condition~(\ref{hausstrat}), since $\partial M$ is a compact $(n-1)$-dimensional manifold. By Example~\ref{croc}, $B_N=\emptyset$ or is a point, so that $M$ satisfies Condition~(\ref{bettistrat0}). By Proposition~\ref{bobord}, $M$ satisfies condition~(\ref{zetastrat}). We can thus apply Theorem~\ref{sph03}. 
	\end{preuve}

\begin{preuve}[ of Theorem~\ref{corobetti}]
	By Corollary~\ref{coq},  for any Morse function $f : M\to \R$, 
	for any non-critical $u\in \R$, 
	$$
	 \forall i\in \{1, \cdots, n\}, \ 
b_i(S_u(M,f))\leq \gamma_i C_i(S_u).$$ 
Condition~(\ref{bettistrat}) and Lemma~\ref{cracboum} imply that $\gamma_i$ is finite. By Theorem~\ref{critical}, 
this implies that 
	$$\forall u\geq -u_0, \ 0\leq  \mathbb E  b(S_u (M,f))-\mathbb Eb_0(S_u)=  \sum_{i=1}^n \mathbb E  b_i(S_u)\leq 
\sum_{i=1}^n\gamma_i\sum_{j\leq n}|\eta_{j,u}|,$$
where $\eta_{j,u}$ satisfies the bound~(\ref{eta}). Now by Theorem~\ref{sph03},
$$\mathbb E b_0(S_u(M,f))=\frac{1}{\sqrt{2\pi}^{n+1}}v(\partial_n M)|u|^{n-1} e^{-\frac12 u^2}+ \delta_u,$$
where $\delta_u$ satisfies~(\ref{deltab}) and $v(M)$ is defined by~(\ref{ve}). 
Hence, there exists $u_2\in \R$ depending on $M$ and $f$ such that
$$\forall u\leq -u_2,\ 
\sum_{j\leq n-1}\gamma_i|\eta_{j,u}|\leq \mathbb E b_0(S_u) e^{-\frac14\theta u^2}.$$

		\end{preuve}
\subsection{Cone spaces}\label{conespace}

In~\cite{adler}, the results for the Euler characteristic hold for a particular type of Whitney stratified sets, the \emph{cone spaces}.
We recall it its definition in a bit more explicit way than the one given by~\cite[8.3.1]{adler}, which
is a bit stronger than the one given by~\cite[3.10.1]{pflaum}:
\begin{definition}\label{conespacedef}
	Let $1\leq \ell\leq k \in \Nn$ and $\widetilde M$ be a $C^k$ manifold.
	\begin{itemize}
		\item A \emph{cone space} $M\subset \widetilde M$ of class $C^\ell$ and \emph{depth $0$}
		is the 	topological
		sum  of		countably many $C^\ell$
		connected submanifolds (without boundary) of $\widetilde M$ together with the stratification  $Z$, the strata of which are given by
		the union of connected components of equal dimension. 
		\item For $d\in \Nn$, a \emph{cone space} of class $\ell$ and
		\emph{depth $d + 1$}  is a stratified space $M\subset \widetilde M$ 
		such that for all $x\in M$ in its stratum $S$, there exist a connected neighborhood $U$ of $x$ in $\widetilde M$, $N\in \Nn$, a compact cone
		space $L\subset \mathbb S^{N-1}\subset \R^{N}$ of class $\ell$ and depth $d$ and finally a $C^\ell$ diffeomorphism 
		$	\varphi : U \to U\cap S \times \mathbb B^{N'}$ such that $$\varphi (U\cap M)=(U\cap S)\times \Cone(L)),$$ where $\Cone(L)=\{tx\in \R^N, t\in [0,1], x\in L\}.$ We also impose that $\varphi_{|U\cap M}$ sends the strata of $M$ homeomorphically onto the natural stratification
		given by the one of $\text{Cone L}$.
	\end{itemize}
\begin{theorem}~\cite[Theorem 3.10.4]{pflaum} A $C^\ell$ cone space is a $C^\ell$ Whitney stratified space. 
\end{theorem}
	\begin{example} The cusp $M=\{(x,y)\in \R^2, y^2=x^3\}\subset \widetilde M=\R^2$ cannot be a $C^1$ cone set. However it is homeomophic to $\Cone(\{\pm 1\})$ and is a Whitney stratified space. A $C^\ell$ manifold of dimension $n$ is a $C^\ell$ cone space of vanishing depth. If $M$ has a boundary and $x\in \partial M$, $x$ then $M$ is locally the product of $\partial M$ and $[0,1]$ which is $\Cone(\{1\})$, hence is a cone space of depth 1. A neighborhood of a vertex of a square is a cone over a quarter of a circle. Using the two first examples, this proves that the square is a smooth cone space.  Similarly, a affine cubes are smooth cone spaces, see~\cite{adler} for other examples. 
	\end{example}
\end{definition}

\subsection{Locally convex sets}
Manifolds with or without boundary and convex polytopes, in particular cubes, belong to a subfamily of Whitney stratified sets and cone spaces which are called in~\cite{adler} \emph{locally convex stratified sets}. 
For these spaces, Morse theory is a bit more explicit.
%
We need further notations. Let $M\subset \widetilde M$ be any subset of a manifold $\widetilde M$, and $x\in M$. Then, the \emph{support cone} $\mathcal S_x M$ is defined by~\cite[(8.2.1)]{adler}:
$$ \mathcal S_x M = \left\lbrace v\in T_x\widetilde M, \ \exists \epsilon>0, \, |\,  \exists c \in C^1 ([0,\epsilon],\widetilde M)\cap C^0([0,\epsilon], M),\  c'(0)=v\right\rbrace.$$
Roughly speaking $\mathcal S_xM$ is the set of directions pointing inwards $M$ from $x$. 
\begin{lemma}\label{iconic}
	Let $M\subset \widetilde M$ be a $C^1$ cone subspace. Then, 
	for any $x\in M$, there exists an integer $N$, a connected neighborhood $U\ni x$ in $\widetilde M$, a diffeomorphism
	$\varphi : U\to \mathbb B^N$ and a connected neighborhood $0\in V\subset T_x S$ such that 
\begin{equation}\label{panurge}\varphi(U\cap M )= (V\cap T_x S) \times (\mathcal S_x\cap \mathbb B^N).
\end{equation}
	\end{lemma}
\begin{preuve}
	By Definition~\ref{conespacedef}, the conclusion of Lemma~\ref{iconic} holds except that 
	$$ \varphi(U\cap M )= (U\cap  S) \times (\Cone(L)).$$ 
	Since $S$ is locally diffeomorphic to $T_xS$, 
	we can change in the latter $S$ into $T_xS$. 
Moreover, 
$$\mathcal S_x M =T_x S \times \{tz, t\in \R_{\geq 0}, z\in L\}\subset T_x\widetilde M\oplus \R^N,$$ 
so that by definition of $\Cone (L)$, (\ref{panurge}) holds.
	\end{preuve}

Assume that the ambient space $\widetilde M$ is equipped with a metric $g$, and let $M\subset (\widetilde M,g)$ be a stratified set. 
The \emph{normal cone} of $M$ at $x\in M$ is defined by~\cite[(8.2.3)]{adler}:
\begin{equation}\label{nxm}
\mathcal N_x M :=\{w\in T_x\widetilde M, \ \forall v\in \mathcal S_x M, \ \langle v,w\rangle \leq 0\}.
\end{equation}
\begin{example} 
	\begin{itemize}
		\item 
		If $M\subset \widetilde M$ is a submanifold, then $\mathcal S_x M = T_x M$ and $\mathcal N_xM\subset T_x \widetilde M$ is the normal bundle of $M$ at $x$. In particular, if $M$ of vanishing codimension, then $\mathcal N_xM=\{0\}$.
		\item If $M\subset \widetilde M$ is a submanifold with boundary and if $x\in \partial M$, then $\mathcal S_x M$ is the half-space in $T_x M$  delimited by $T_x \partial M$  in the inward direction, and $\mathcal N_xM$ is the convex cone generated by an outward normal vector $n_x$ in $T_xM$ (orthogonal to $T_x\partial M$) and the normal bundle of $M$ at $x$ in $\widetilde M$. In particular, if $M$ is of vanishing codimension, $\mathcal N_xM= \R_{\geq 0} n_x$.
		\item If $M \subset \R^2$ is a rectangle and $x$ is a vertex, then $\mathcal S_x M$ is the the cone of directions parallel to the inner quartant at $x$, and $\mathcal N_xM= -\mathcal S_xM$.
	\end{itemize}
\end{example}

\begin{definition}~\cite[Definition 8.2.1]{adler}\label{locconv}
	A Whitney stratified set $M$ is \emph{locally convex} if 
	for any $x\in M$, $\mathcal S_x M$ is convex. 
\end{definition}
\begin{example}
	Manifolds with or without boundary and affine convex polytope are locally convex. Note that a plane polytope with a concave angle is not locally convex at the concave summit, but is a smooth cone space. 
\end{example}


The following lemma generalizes Example~\ref{croc} in this particular subfamily of stratified sets.
\begin{lemma}\label{Gogo}
	Let $M\subset (\widetilde M,g)$ be a locally convex $C^1$ cone space of dimension $n$, $x\in S\subset M$, and $\nu \in T^\perp_x M$. Then,
	$$\left\lbrace
	\begin{array}{llll}
	(A_N,B_N)(x,\nu)&\sim_{hom}& (\mathbb B^{n-\dim S}, \emptyset) &\text{ if } 
	-\nu \in \mathcal N_x M \\
	 (A_N,B_N)(x,\nu)&\sim_{hom}& (\mathbb B^{n-\dim S}, \mathbb B^{n-\dim S-1})  &
	\text{ if }  -\nu\notin \mathcal N_x M.
	\end{array}\right.
	$$
	If $p$ is a critical point for $f$ with $f(p)=u$, then
		$$\left\lbrace
	\begin{array}{llll}
	(A,B)&\sim_{hom}& (\mathbb B^{n}, \mathbb B^{n-\dim S}\times B_T) &\text{ if } 
	-\nabla f(p)\in \mathcal N_x M \\
	(A,B)&\sim_{hom}& \left(\mathbb B^{n},( \mathbb B^{n-\dim S}\times B_T) \cup (\mathbb B^{n-\dim S-1}\times A_T)\right)  & \\
&	&\text{ and }
	b_0(S_{u+\epsilon})=b_0( S_{u-\epsilon}) &
	\text{ if }  -\nabla f(p)\notin \mathcal N_x M.
	\end{array}\right.
	$$
		In particular, if $p$ has vanishing (tangent) index, then
				$$\left\lbrace
		\begin{array}{lllll}
		(A,B)&\sim_{hom}& (\mathbb B^{n}, \emptyset)&\text{ and }
		S_{u+\epsilon}\sim_{hom} S_{u-\epsilon}\sqcup \mathbb B^n 
		 &\text{ if } 
		-\nabla f(p)\in \mathcal N_x M \\
		(A,B)&\sim_{hom}& (\mathbb B^{n}, \mathbb B^{n-1})  &\text{ and }
		b_0(S_{u+\epsilon})=b_0( S_{u-\epsilon}) &
		\text{ if }  -\nabla f(p)\notin \mathcal N_x M.
		\end{array}\right.
		$$
	\end{lemma}
In Figure~\ref{tore-bord}, $q$ is an example of the last situation, and $p$ for the penultimate situation. 
\begin{preuve} By Lemma~\ref{iconic}, we can assume that 
	$M=\R^{\dim S}\times \mathcal S_x\subset \R^{\dim S+N}$ equipped with the standard scalar product  in the neighborhood of $x$.  Let $f_\nu(y)=\langle y,\nu\rangle $.  Then, $df_\nu (0) = \nu^*$, and $f^{-1}[-\epsilon, \epsilon]$ is a linear band, hence convex, of vanishing codimension containing $0$ in its interior. Hence, since $\mathcal S_x$ has interior being of dimension $n-\dim S$ and is convex, as any transverse ball $D_x$ (see paragraph~\ref{dix}), 
	so that 
	$$f_\nu^{-1}[-\epsilon, \epsilon]\cap N_x\cap B(0,\delta)$$ is a convex subset of $\R^{\dim S+N}$ with non-empty interior containing a ball of dimension $n-\dim S$, so that it is homeomorphic to a ball of dimension $n-\dim S$. Hence, for any $\nu$, $$A_N(x,\nu)\sim_{hom} \mathbb B^{n-\dim S}.$$
	
	Assume now that $\nu \in -\mathcal N_x M$. Then, $\langle \nu,v\rangle \geq 0$ for any $v\in \mathcal S_x M$. This implies that $\{y\in \R^{\dim S+N}, \langle y,\nu\rangle <0\}$ is a half space whose intersection with $\mathcal S_x M$ is empty. In particular, $B_N(x,\nu)=\emptyset$. 
	
Assume next that $\nu \notin -\mathcal N_xM$. Then, there exists $v\in \mathcal S_xM$, such that  $\langle \nu,v\rangle < 0$. Hence, for any $\epsilon>0$, the affine hyperplane 
	$\{y\in \R^{\dim S+N}, \langle y,\nu\rangle =-\epsilon\}$ intersects $\mathcal S_x M$ in its interior. By the same arguments given for $A_N$, this implies that $$f_\nu^{-1}(-\epsilon)\cap N_x\cap B(0,\delta)$$ 
is homeomorphic to a ball of dimension $n-\dim S-1$, hence the same for $B_N(x,\nu)$.

The two first general assertion concerning $(A,B)$ are now a direct consequence of its definition, as the first of the last pair of assertions. The last assertion is due to the fact that the connected handle $A$ is attached through $B$ which is connected.
\end{preuve}
The following corollary is a  generalization of
Corollary~\ref{coromilnor}.
\begin{corollary}\label{comparaison} Let $M\subset (\widetilde M,g)$ be a locally convex $C^2$ cone space of dimension $n$. Then, for any Morse function $f : M\to \R$ and any real $u$ which is not a critical value of $f$,
	$$ 0\leq 
		\sum_{x\in \crit_0(S_u(M,f))}
	{\bf 1}_{\{-\nabla f(x) \in \mathcal N_x M\}} - 
	N(S_u(M,f))\leq \sum_{i\geq 1} C_i(S_u(M,f)).
	$$
\end{corollary}
\begin{remark}\label{tark}
	\begin{enumerate}
		\item 
	If $M$ is a manifold with boundary, F. Laudenbach explained us how to use use~\cite{laudenbach} to prove that at a critical point $p$ on the boundary with $\nabla f(x)$ in the direction of $n_x$, then   $S_{f(p)+\epsilon}\sim_{homeo} S_{f(p)-\epsilon}$, so that Corollary~\ref{comparaison} should hold for $N_{\mathbb B}$ (more correctly, an homeomorphic version of it) instead of $N(S_u)$. It is very likely that the same holds for general locally convex cone sets. 
	\item Corollary~\ref{comparaison} is not true for $N_{\mathbb S^{n-1}}(Z_u) $ or $N(Z_u)$ for stratified sets which are not compact manifolds without boundaries. Indeed in the example of Figure~\ref{tore-bizarre}, passing the critical point $q$, where the gradient points outwards, changes $Z_{d-\epsilon}\sim \mathbb S^1$ but $Z_{d+\epsilon}\sim[0,1]$. Moreover $b_0(Z_u)$ jumps from 1 to 0 after the highest critical point. 
	\end{enumerate}
\end{remark}
\begin{preuve} [ of Corollary~\ref{comparaison}]
	By Lemma~\ref{Gogo}, any local minimum in its stratum creates a connected component of $S_u$
	if $-\nabla f (x)\in \mathcal N_x M$, and in this case the component is homeomorphic to a ball of maximal dimension, and in the other case, 
	the number of components of the upper level is the same as the lower level. 
	 Moreover, any critical point of positive index cannot create a component, since $B\not=\emptyset$ in this case. 
	Hence, 
		$$ 	
		 N(S_u)\leq 
			\sum_{x\in \crit_0(S_u)}
	{\bf 1}_{\{-\nabla f(x) \in \mathcal N_x M\}}.
	$$
By Lemma~\ref{relative}, a critical point can kill at most one connected component. Hence, 
$$N(S_u)\geq 
	\sum_{x\in \crit_0( S_u)}
{\bf 1}_{\{-\nabla f \in \mathcal N_x M\}} - 
\sum_{i\geq 1} C_i(S_u).$$
These two pairs of inequalities prove the result.
\end{preuve}

%

\subsection{The refinement }

In this paragraph we want to prove the following quantitative version of Theorem~\ref{sph2z}, which  is a 
improvement of Theorem~\ref{sph03}. On the contrary to the latter, Theorem~\ref{sph2} uses the main result of~\cite{adler}, holds only for locally convex cone spaces. 
\begin{theorem}\label{sph2}	Let $\widetilde M$ be a $C^3$ manifold of dimension $n\geq 1$,  $M\subset \widetilde M$ be a compact locally convex $C^2$  cone space of dimension $n$ satisfying condition~(\ref{grass}) (very gentle boundaries),
	$ \widetilde f: \widetilde M \to \R$ be a random centered Gaussian field satisfying conditions (\ref{C1}) (regularity), (\ref{C2}) (non-degeneracity) and (\ref{C4}) (constant variance) defined below, $f=\widetilde f_{|M}$ and $g$ be the metric induced by $f$ and defined by~(\ref{metric}). 
Then, 
	$$\forall u\geq u_1,\ \mathbb E N(E_u(M,f)) =
	\sum^n_{i=0} \frac{1}{\sqrt{2\pi}^{i+1}}\mathcal L_i H_{i-1} (u)e^{-\frac{u^2}2},
	+r_u,$$
	where the constants $(\mathcal L_k)_k$ and the Hermite polynomials $(H_k)_k$ are defined below by~(\ref{killing}) and~(\ref{hermite}), and where 
	\begin{equation}\label{rume}
	\forall u\geq u_1,\ |r_u|\leq 
	\frac{1}{\rho}Q_N(\sigma^{1/2}, \sigma^{-1/2}, s)	\sup_{0\leq j\leq n} v(\partial_j M)u^{N^2}e^{-\frac{1}2u^2(1+\theta)}.
	\end{equation}
Here $v(M)$, $u_1, \sigma, \rho,\theta $ are defined by~(\ref{ve}), (\ref{uun}), ~(\ref{sigma}), ~(\ref{rho}), (\ref{cem}), ~(\ref{psi}) and (\ref{phi}) and $Q_N$ is a polynomial depending only on $N$  with non-negative coefficients.

 If $M$ is a compact $C^2$ manifold without boundary,	the same holds for 
	$N_{\mathbb B^n}(E_u(M,f))$, 	$N_{\mathbb S^{n-1}}(Z_u)$ and $N(Z_u)$ instead of 
	$ N(E_u),$ after changing the polynomial $Q_N$.
\end{theorem} 

Theorem~\ref{sph2} is a consequence of the following Theorem~\ref{sph1} and the main result of~\cite{adler}, namely Theorem~\ref{euler} below. 
\begin{theorem}\label{sph1}
	Under the hypotheses of Theorem~\ref{sph2},
	then 
	$$\forall u\geq u_0, \ \mathbb E N(E_u(M,f)) =
	\mathbb E \chi (E_u)+r_u$$
where 
	\begin{equation}\label{rume2}
\forall u\geq u_0, \	|r_u|\leq 
\frac{1}{\rho}Q_N(\sigma^{1/2}, \sigma^{-1/2}, s)
	\sup_{0\leq j\leq n}v(\partial_j M)|u|^{N^2}e^{-\frac{1}2u^2(1+\theta)}.
	\end{equation}
	Here, $u_0$ and $s$ are defined by~(\ref{phi}), and $\rho, \sigma^{-1}$ and $\theta$ by~(\ref{psi}).

	If $M$ is a compact manifold without boundary, the same holds for 
	$N_{\mathbb B^n}(E_u)$, 	$N_{\mathbb S^{n-1}}(Z_u)$ and $N(Z_u)$ instead of 
	$ N(E_u),$ after changing the polynomial $Q_N$.
\end{theorem} 
\begin{lemma}~\cite[Corollary 9.3.3]{adler}\label{fofo} Under the hypotheses of Theorem~\ref{sph2}, 
for any Morse function $f: M\to \R$ and any non-critical $u\in \R$,
	$$ \chi(S_u(M,f))= \sum_{j=0}^n 
	\sum_{x\in \crit (S_u(\partial_j M,f))}
	(-1)^{\Ind(x)} {\bf 1}_{\{\nabla f (x) \in -\mathcal N_x M\}},$$
	where $\mathcal N_x$ is defined by~(\ref{nxm}).
	\end{lemma}
	This equality is the locally convex stratified version of~(\ref{chi}).
\begin{preuve}[ of Theorem~\ref{sph1}]
 By Corollary~\ref{comparaison},
		$$ 0\leq \sum_{j=0}^n 
	\sum_{x\in \crit_0(\partial_j M,f)}
	{\bf 1}_{\{-\nabla f \in \mathcal N_x M\}} - 
	N(S_u(M,f))\leq \sum_{i\geq 1} C_i(S_u(M,f)),
	$$
	so that by Lemma~\ref{fofo}
	$$ \big|\chi(S_u)-
	N(S_u(M,f))\big|\leq 2\sum_{i\geq 1} C_i(S_u(M,f)).
	$$
	By Corollary~\ref{critical}, for $u\leq -u_0,$ the right-hand side is bounded by
	$$ 	
\frac{2}{\rho}Q_N(\sigma^{1/2}, \sigma^{-1/2}, s)	\sum_{j=0}^nv(\partial_j M  ) |u|^{N^2}e^{-\frac{1}2u^2(1+\theta)},
	$$
	where $\rho, \sigma^{-1}$ and $\theta$ are given by~(\ref{phi}) and $Q_N$ is a polynomial depending only on $N$. Hence, the result.

Assume now that $M$ is a manifold without boundary. Then, for any critical point of $f$, since $\mathcal N_x$ is the normal bundle at $x$, $-\nabla f(x)\in \mathcal N_x$ so that 
$$
{\bf 1}_{\{x\in \crit(M,f), \ -\nabla f \in \mathcal N_x M\}} = {\bf 1}_{x\in \crit(M,f)}.$$
Moreover by Corollary~\ref{coromilnor},
$$ 0\leq 
C_0(M)
- 
N_{\mathbb B^n}(S_u(M,f))\leq \sum_{i\geq 1} C_i(S_u).
$$
The sequel is the same as above in the general case. Corollary~\ref{pasteque} provides the analogous argument for $Z_u$.
\end{preuve}
Theorem~\ref{sph1} must be associated to the following Theorem~\ref{euler}, which is a exact formula computing the average Euler characteristic of $E_u(M,f)$ in this context of a a regular locally convex cone space:
\begin{theorem}\label{euler}\cite[Theorem 12.4.2 and Remark 12.4.3]{adler}
Under the hypotheses of Theorem~\ref{sph2}, 
	$$\forall u\in \R, \ \mathbb E \chi(E_u(M,f))= 
	\sum^n_{j=0} \frac{1}{\sqrt{2\pi}^{j+1}} \mathcal L_j H_{j-1} (u)e^{-\frac{u^2}2}.$$
\end{theorem}
Here, for every $ k\in \{0, \cdots, n\},$ the \emph{Lipschitz-Killing curvature} $\mathcal L_k$  is defined by
\begin{align}\label{killing}
	\mathcal L_k (M) =&\sum_{j=k}^n \frac1{\sqrt{2\pi}^{j-k}}
	\sum_{\ell=0}^{\lfloor(j-k)/2\rfloor} 
	C_{n-j, j-k-2\ell}
	\frac{(-1)^\ell}{\ell! (j-k-2\ell)!} \\ \nonumber
	& \int_{x\in \partial_j M}
	\int_{\nu_{n-j}\in S(T_x\partial_j M^\perp)}{\bf 1}_{\{\nu_{n-j }\in -\mathcal N_x\}}
	\Tr^{T_x \partial_jM} \left(R^\ell S^{j-k-2\ell}_{\nu_{n-j}}\right)
	 d\vol_{n-j-1}(\nu_{n-j}) d\vol_g
	(x).
\end{align}	
Let us explain the notations of Theorem~\ref{euler}. First, 
$H_j$ is the $j$th Hermite polynomial, that is:
\begin{equation}\label{hermite}
 \forall x\in \R, \ H_{-1}(x) = \sqrt{2\pi}\Psi(x) e^{x^2/2} \text{ and }\forall j\geq 0, \ H_j (x) = (-1)^j e^{x^2/2} \frac{d^j}{dx^j} (e^{-x^2/2}),
\end{equation}
where $\Psi(x) = \frac{1}{\sqrt{2\pi}}\int^{\infty}_x e^{-\frac12 u^2}du.$
Note that 
$$ \forall j\geq 1, \ H_j(x) = j!\sum_{\ell=0}^{\lfloor j/2\rfloor}
\frac{(-1)^\ell x^{j-2\ell}}{\ell! (j-2\ell)! 2^\ell},$$
so that 
$
\forall j\geq 0, H_{j}(u)\underset{u\to +\infty}{\sim}
u^{j}.
$
Moreover, 
$$\forall m,i \in \Nn, \, C_{m,i} = \left\lbrace\begin{array}{ccc} &\frac{(2\pi)^{i/2}}{s_{m+i}} &\text{ if } m>0\\
 &1& \text{ if } m=0 
 \end{array}\right., $$
 where $\displaystyle s_m = \frac{2\pi^{\frac{m}2}}{\Gamma(\frac{m}2)}=\vol_{g_0} \mathbb S^{m-1}.$

\begin{preuve}[ of Theorem~\ref{sph2}] This is an immediate consequence of Theorems~\ref{sph1} and~\ref{euler}. 
\end{preuve}
\begin{example}\label{exkilling}
		Let $f : \R^n\to \R$ be a centered Gaussian field satisfying conditions (\ref{C1})(regularity), (\ref{C2})(non-degeneraticity) and (\ref{C3}) (stationarity),
		and $M$ be a compact open set with smooth boundary $\partial M$.
		In this case, $R=0$ 
		and for any $x\in \partial M$, ${\bf 1}_{\mathcal N_x}={\bf 1}_{\nu_1=n_x}$, where $n_x$ denotes the outward unit normal vector to $\partial M$. Then,
			\begin{eqnarray*}
				 \mathcal L_n (M) &=&
				\vol_g{M} \\
\text{ and }			\forall k\in \{0, \cdots, n-1\}, \
			\mathcal L_k (M) &=&
\frac{1}{s_{n-k}}
			\frac{1}{(n-1-k)!} 
			 \int_{\partial M}
			\Tr^{T_x \partial M} \left(S^{n-1-k}_{-n_{x}}\right)		d\vol_g (x). 
		\end{eqnarray*}
	\end{example}

\subsection{The asymptotic of  $c_{Z}(u)$ }

Fix $n\in \Nn^*$.  In the sequel, $\forall r>0, \ B_r=r\mathbb B^n\subset \R^n.$
 We begin by recall the main result of~\cite{nazarov2}.
%
\begin{theorem}\cite[Theorem 1.1]{nazarov2}\label{nst} Let $f : \R^n \to \R $ be a random centered Gaussian field satisfying conditions~(\ref{C1}) (regularity), (\ref{C2}) (non-degeneraticity), (\ref{C3}) (stationarity) and~(\ref{C5}) (ergodicity).
	Then, there exists a non negative constant $c_{Z}(u)$ such that 
	$$ \frac{\mathbb E N(Z_u(B_r, f))}{\vol_{g_0} (B_r)}\underset{r\to +\infty}{\to} c_{Z}(u).$$
\end{theorem}
The conditions for this theorem are in fact milder, see~\cite{nazarov2}.
The proof of Theorem~\ref{ns2} is not a direct consequence of Theorem~\ref{sph1}. Indeed, 
the latter holds for $E_u$ but not for $Z_u$. 
For the proof of Theorem~\ref{ns2}, we will need the following simple lemma. 
\begin{lemma}\label{rescale}
	Let $M\subset \R^n$ be a $C^2$ compact
	codimension 0 submanifold with $C^2$ boundary $\partial M.$ 
	Assume that $\R^n$ is equipped with a stationary metric $g$. 
	For any $r>0$, let $M_r=rM$ and
	\begin{equation}\label{esser}
	S_r: T\partial M_r \times T\partial M_r \to T^\perp \partial M_r
	\end{equation}
	 be the second fundamental form  associated to the pair $(\partial M_r,\R^n)$, defined by~(\ref{second}), where the affine space is equipped with the metric $g$. Then,
	$$ \forall r>0, x\in M, \ S_r(rx)= \frac{1}r S_1(x),$$
	where we identify $T_{rx}\partial M_r $ with $ T_x \partial M$.
	In particular, for all $r>0$, 
	$ s(M_r)=\frac{1}rs(M), $ where
	$s$	is defined by~(\ref{esse}).
\end{lemma}
\begin{preuve} Since $g$ a constant metric over $\R^n$, 
	$\nabla = d$. 
	Let  $x\in \partial M$
	and $X,Y : \partial M \to \R^n$ two tangent vector fields of $T\partial M$ near $x$. Fix $r>0$. Then, $X_r(\cdot):=X(\cdot/r )$ and $Y_r(\cdot):= Y_r(\cdot/r)$) are tangent vector fields of $\partial M_r$ near $rx$, and $$\nabla_{X}Y (x)= rd_{X_r}Y_r(x)=r\nabla_{X_r} Y_r (rx),$$
	so that if $P_{T^\perp\partial_r M} : \R^n \to T^\perp\partial M_r$ denotes the orthogonal (for $g$) projection onto the normal bundle of $\partial  M_r$, then 
	$$  S_1(X,Y)=P_{T^\perp M}\nabla_{X} Y (x) =r P_{T^\perp\partial_r M}\nabla  _{X_r }Y_r(rx)= rS_r(X_r,Y_r)=rS_r(X,Y),$$
	where we identified $X$ and $Y$  with $X_r$ and $Y_r$ as vectors. Hence, the result.
\end{preuve}
\begin{preuve}[ of Theorem~\ref{ns2}] Recall that $B_r=r\mathbb B^n \subset \R^n$. 
	By Lemma~\ref{rescale}, $$\forall r>0, \ s_r= \frac1r s_1.$$ 
	Moreover since $R=0$, $\rho$
	defined by~(\ref{rho}) is equal to 1, and $\sigma^{-1}$ defined by~(\ref{sigma}) equal to 1 as well. In particular, 
	\begin{align}\label{uzb}
	\forall r\geq 1,\ 
	u_0(B_r) = (1+s_1/r)\leq 1+s_1,
	\end{align}
	where $u_0$ is defined by~(\ref{uzero}),
	and $$\theta(B_r) = \frac{1}{(1+s_1/r)^2}\geq \frac{1}{(1+s_1)^2},$$ 
	where $\theta$ is defined by~(\ref{cem}).
Hence, $$u_1(B_r)\leq (n^2+2)(1+s_1)^2,$$
where $u_1$ is defined by~(\ref{uun}).
	Since $g=d^2 e(0)$, 
	$$d\vol_g = (\det d^2 e (0))^{1/2} d\vol_{g_0},$$ where the determinant is computed in the standard basis of $\R^n$.
	Recall that $\overline{B_r}$ is stratified as $B_r \cup \partial B_r$.
By Corollary~\ref{dauphin}, 
	\begin{equation}\label{coquette}
		 \forall u\in \R, \ |b_0(Z_u(\overline{B_r}))- C_0(S_u(B_r))|\leq 
	\zeta_0(C(S_u(\overline{B_r}))-C_0(S_u(B_r)),
	\end{equation}
	where $\zeta_0$ is defined by~(\ref{zeta}). By Proposition~\ref{bobord}, $\zeta_0$ is bounded by a constant depending only on $n$. 	Moreover by Lemma~\ref{fofo}
	\begin{equation}\label{cocorico}
	\forall u\in \R, \ |\chi(S_u(\overline{B_r}))-C_0(S_u({B_r}))|
	\leq + C(S_u(\partial B_r))+ \sum_{i=1}^n C_i(S_u(B_r)).
	\end{equation}
	Using Theorem~\ref{euler} and Example~\ref{exkilling},  
	(\ref{coquette}) and (\ref{cocorico}) imply that 
	for any 	$u\in \R$,
	\begin{align}\label{cocorico}
\frac{\mathbb E N(Z_u(\overline{B_r}))}{\vol_{g} B_r}& =\frac{1}{\sqrt{2\pi}^{n+1}} H_{n-1}(u)e^{-\frac12 u^2}+e^{-\frac12 u^2}
\\
&
\sum_{k=0}^{n-1} \frac{1}{\sqrt{2\pi}^{k+1}} H_{k-1}(u)\frac{1}{s_{n-k}}
\frac{1}{(n-1-k)!} 
\int_{\partial B_r}
\Tr^{T_x \partial B_r} \left(S^{n-1-k}_{r|-n_{x}}\right)		\frac{d\vol_g (x)}{\vol_{g} B_r}+\mu_u,
	\end{align} 
where $S_r$ is defined by~(\ref{esser}) and 
\begin{align}\label{mumumu}
 |\mu_u|\leq \frac{1+\zeta_0}{\vol_{g} B_r} \left(\mathbb E C(\partial B_r) +\sum_{i=1}^n \mathbb E C_i(B_r)\right).
 \end{align}
  By
Lemma~\ref{rescale}, for any $r\geq 1$, 
$$\Tr^{T_x \partial B_r} (S^{n-1-k}_{r|-n_{x}})= \Tr^{T_x \partial B_r} \frac{1}{r^{n-1-k}}(S^{n-1-k}_{1|-n_{x}}).$$
Hence, taking $r\to +\infty$ in~(\ref{cocorico}) kills the sum of boundary terms in the above equation. 
By Theorem~\ref{critical} applied to $B_r$ and Proposition~\ref{critz} applied to $\partial B_r$,
for $r\geq 1$, (\ref{mumumu}) gives, using~(\ref{uzb}), 
\begin{align*}\forall u\leq -(1+s_1), \ 
\frac{|\mu_u|}{1+\zeta_0}&\leq  	
Q_n(1,1,0)|u|^{n^2}e^{-\frac{1}2u^2(1+\theta)}+\\
&\frac{2\vol_g (\partial B_r  )}{\vol_g B_r}|u|^{n-3} e^{-\frac12 u^2}\left(
	\frac{|u|}{\sqrt{2\pi}^{n}}+ P_n(1, 1,s_1)
	\right).
\end{align*}
When $r$ goes to $+\infty$, the second term vanishes, hence the result.
\end{preuve}

\bibliographystyle{amsplain}
\bibliography{excursion.bib}

\noindent Univ. Grenoble Alpes, Institut Fourier \\
F-38000 Grenoble, France \\
CNRS UMR 5208  \\
CNRS, IF, F-38000 Grenoble, France

\end{document}